# New Formulas and Methods for Interpolation, Numerical Differentiation and Numerical Integration

**M. Ramesh Kumar**
Phone No: +91 9840913580
Email address: ramjan_80@yahoo.com
Home page url: http//ramjan07.page.tl/

**Abstract**

We present a new formula for divided difference and few new schemes of divided difference tables in this paper. Through this, we derive new interpolation, numerical differentiation and numerical integration formulas with arbitrary order of accuracy for evenly and unevenly spaced data. First, we study the new interpolation formula, which generalizes both Newton's and Lagrange's interpolation formula with the new divided difference table for unevenly spaced points and using this; we derive other interpolation formulas, in terms of differences and divided differences for evenly spaced data. Second, we study two new different methods of numerical differentiation for both evenly and unevenly spaced points without differentiating the interpolating polynomials or the use of operators. Third, we derive new numerical integration formulas using new differentiation formulas and Taylor formula for both evenly and unevenly spaced data. Basic computer algorithms for few new formulas are given. In comparison to former polynomial interpolation, numerical differentiation and numerical integration formulas, these new formulas have some new featured advantages for approximating functional values, numerical derivatives of higher order and approximate integral values for evenly and unevenly spaced data.

## 1. Introduction

Polynomial Interpolation theory has a number of important uses. Its primary use is to furnish some mathematical tools that are used in developing methods in the areas of approximation theory, numerical integration and the numerical solution of differential equations [1]. The general problem of interpolation consists, then, in representing a function, known or unknown, in a form chosen in advance, with the aid of given values which this function takes for definite values of the independent variable [15]. A number of different methods have been developed to construct useful interpolation formulas for evenly and unevenly spaced points. Newton's divided difference formula [1,5] and Lagrange's formula [3,6,14] are the most popular interpolation formulas for polynomial interpolation to any arbitrary degree with finite number of points.

Lagrange interpolation is a well known, classical technique for interpolation. Using this; one can generate a single polynomial expression which passes through every point given. This requires no additional information about the points. This can be really bad in some cases, as for large numbers of points we get very high degree polynomials which tend to oscillate violently, especially if the points are not so close together. It can be rewritten in two more computationally attractive forms: a modified Lagrange form [2] and a bary centric form [2, 10]. Newton's formula for constructing the interpolation polynomial makes the use of divided differences through Newton's divided difference table for unevenly spaced data [1]. Based on this formula, there exists many number of interpolation formulas using differences through difference table, for evenly spaced data. The best formula is chosen by speed of convergence, but each formula converges faster than other under certain situations, no other formula is preferable in all cases. For example, if the interpolated value is closer to the center of the table then we go for any one of central difference formulas, (Gauss's, Stirling's and Bessel's etc) depending on the value of argument position from the center of the table. However, Newton interpolation formula is easier for hand computation but Lagrange interpolation formula is easier when it comes to computer programming. In this paper, we study a new interpolation formula which generalizes both Newton's and Lagrange's interpolation formula and its various forms with new scheme of divided difference tables for evenly and unevenly spaced data. Also, we study the comparisons of new interpolation formulas with the former interpolation formulas based on differences.

Numerical approximations to derivatives are used mainly in two ways. First, we are interested in calculating derivatives of given data that are often obtained empirically. Second, numerical differentiation formulae are used in deriving numerical methods for solving ordinary and partial differential equations [1]. A number of different methods have been developed to construct useful formulas for numerical derivatives. Most popular of the techniques are finite difference type [4,14], polynomial interpolation type [3,4,12], method of undetermined coefficients [9], and Richardson extrapolation [3,4,11,14,16]. The range of numerical differentiation formulae based on different interpolating polynomials and finite difference is limited, as a rule, to finite points of interpolation. All available formulae known at the present moment are obtained for a certain concrete limited number of interpolation points. It can be explained by the fact that the procedure of the finding of the interpolating polynomial coefficients in the case of the arbitrary number of interpolation points is quite uncomfortable and requires difficult calculations.

In Ref [13], the explicit finite difference formulas on numerical differentiation for equally and unequally spaced data are presented. The formulas, given in this paper are very useful to find higher order derivatives directly, but the calculation burdens of the coefficients of explicit differences in the formulas increases when the number of data increases. In Ref [7,8], central difference formulas for first and second order derivatives for finite and infinite number of data are discussed. The coefficients of differences in the formula are easy to implement and take less computation burden. They are calculated by Vandermonde determinant through Crammer's rule. But, this method of undetermined coefficients is limited to first and second order formula in this paper. Thus, the higher order differentiation formulas using central differences are still to be developed for arbitrary number of data. Also, the general formula for numerical differentiation for evenly and unevenly spaced data with less computational burden is still questioned.

There exists number of methods to solve the problem of numerical integration. Few popular methods are Simpson rule [1,5], Newton cotes integration [1,6,9,16], Romberg integration [1,3,5,16], Gaussian quadrature [1,3,5] and Automatic integration [1]. All available formulas are presented only to evenly spaced data and limited to certain order of accuracy. The problem of numerical integration for unevenly spaced data is still questioned. Also, the numerical integration formula for evenly spaced data to any arbitrary order of accuracy is still to be developed.

In Section 2, a new formula for divided differences in terms of polynomial is derived and new schemes of divided difference tables are introduced for evenly and unevenly spaced data. In Section 3, the new interpolation formulas are presented in various forms for both evenly and unevenly spaced. Further, comparisons of the new interpolation formulas with former interpolation formulas (Newton's, Gauss, Stirling, etc) are discussed. In Section 4, two new methods of higher order numerical differentiation formulas are studied for evenly and unevenly with various remainder terms. In Section 5, the direct formulas for numerical differentiation are given in terms of functional values and these new formulas are compared with former differentiation formulas based on interpolating polynomial and Theorems in Ref [13]. In Section 6, the new numerical integration formulas are derived using new differentiation formulas, presented in section 5 and Taylor formula. In Section 7, is devoted to a brief conclusion.

## 2. Formula for divided difference

**Definition 2.1**

The $r^{\text{th}}$ divided difference of polynomial function $P(x)$ at the points $x, x_0, x_1, \ldots x_{r-1}$ is a polynomial in $x$, so we call it as divided difference polynomial of order r of $P(x)$. It is denoted by $P[x, x_0, x_1, \ldots x_{r-1}]$.

**Theorem 2.1.** *Let $x_0, x_1, \ldots x_{r-1}$ are 'r' numbers and $x_{r,} x_{r+1}, \ldots, x_n$ are $(n-r+1)$ distinct numbers in the interval $[p,q]$, $r \le n$, $x \in [p,q]$ and $f \in C^{n+1}[p,q]$, then*

$$f[x, x_0, x_1, \ldots, x_{r-1}] = \sum_{i=r}^{n} f[x_i, x_0, x_1, \ldots, x_{r-1}] \prod_{\substack{j=r \\ i \ne j}}^{n} \frac{x - x_j}{x_i - x_j} + f[x, x_0, x_1, \ldots, x_n] \prod_{i=r}^{n} (x - x_i) \qquad (2.1)$$

**Proof.**

Let $P_n(x)$ is a polynomial of degree $\le n$ in '$x$' that approximates the function $f$ and takes the functional values $f(x_0), f(x_1), \ldots f(x_n)$ for the arguments $x_0, x_1, \ldots x_n$ respectively. It can be written as

$$P_n(x) = a_0 + a_1 x + a_2 x^2 + \ldots + a_n x^n, \qquad \text{All } a's \in R \qquad (2.2)$$

Then, we can write $r^{th}$ order divided difference of $P_n(x)$ at the points $x_0, x_1, \ldots x_{r-1}$ in terms of '$x$' as in the following form

$$P_n[x, x_0, x_1, \ldots x_{r-1}] = \hat{a}_0 + \hat{a}_1 x + \hat{a}_2 x^2 + \ldots + \hat{a}_{n-r} x^{n-r} = \tilde{P}_{n-r}(x) \text{ (Say), All } \hat{a}'s \in R \quad (2.3)$$

Now $\tilde{P}_{n-r}(x)$ is a polynomial of degree $\leq n-r$. To interpolate it, we can use remaining points $x_{r,} x_{r+1}, \ldots, x_n$. Also, we know that, for polynomial of degree $\leq n-r$

$$\tilde{P}_{n-r}[x, x_r, x_{r+1}, x_{r+2}, \ldots, x_n] = 0 \quad (2.4)$$

So we can express,

$$\tilde{P}_{n-r}(x) = \sum_{i=r}^{n} A_i \prod_{\substack{j=r \\ i \neq j}}^{n} (x - x_j) \quad (2.5)$$

Using (2.3) in (2.5)

$$P_n[x, x_0, x_1, \ldots x_{r-1}] = \sum_{i=r}^{n} A_i \prod_{\substack{j=r \\ i \neq j}}^{n} (x - x_j) \quad (2.6)$$

At $(x_k, y_k)$, where $k = r, r+1, r+2, \ldots, n$

$$A_k = \frac{P_n[x_k, x_0, x_1, \ldots x_{r-1}]}{\prod_{j=r, j \neq k}^{n} (x_k - x_j)}, \quad k = r, r+1, r+2, \ldots, n \quad (2.7)$$

Substituting (2.7) in (2.6)

$$P_n[x, x_0, x_1, \ldots, x_{r-1}] = \sum_{i=r}^{n} P_n[x_i, x_0, x_1, \ldots, x_{r-1}] \prod_{\substack{j=r \\ i \neq j}}^{n} \frac{x - x_j}{x_i - x_j} \quad (2.8)$$

The above equation is polynomial form of $r^{th}$ divided difference of a polynomial. Since $P_n(x)$ approximates the function $f(x)$ on [$p$, $q$], so we have the following two equations (2.9) and (2.10).

$$f(x) = P_n(x) + f[x, x_0, x_1, \ldots, x_n] \prod_{i=0}^{n} (x - x_i) \quad (2.9)$$

$$f[x_i, x_0, x_1, x_2 \ldots, x_{r-1}] = P_n[x_i, x_0, x_1, x_2 \ldots x_{r-1}], \quad i = r, r+1, r+2, \ldots, n \quad (2.10)$$

Using (2.9), we can write

$$f[x, x_0, x_1, \ldots, x_{r-1}] = P_n[x, x_0, x_1, \ldots, x_{r-1}] + f[x, x_0, x_1, \ldots, x_n] \prod_{i=r}^{n} (x - x_i) \quad (2.11)$$

Using the equation (2.8) in (2.11), we get

$$f[x, x_0, x_1, \dots, x_{r-1}] = \sum_{i=r}^{n} P_n[x_i, x_0, x_1, \dots, x_{r-1}] \prod_{\substack{j=r \\ i \neq j}}^{n} \frac{x - x_j}{x_i - x_j} + f[x, x_0, x_1, \dots, x_n] \prod_{i=r}^{n} (x - x_i) \tag{2.12}$$

Using equation (2.10) in (2.12), we obtain equation (2.1).

**Note 2.1.** As a notation, $\{a, b, c, \dots\}$ denotes the smallest interval containing all of real numbers $a, b, c, \dots$ [1].

**Corollary 2.1.1.** *Let* $x_0, x_1, \dots x_{r-1}$ *are 'r' numbers and* $x_r, x_{r+1}, \dots, x_n$ *are* $(n - r + 1)$ *distinct numbers in the interval* $[p, q]$, $r \leq n$, $x \in [p, q]$ *and* $f \in C^{n+1}[p, q]$, *then for some* $\xi \in \{x, x_0, x_1, \dots x_n\}$ *with*

$$f[x, x_0, x_1, \dots, x_{r-1}] = \sum_{i=r}^{n} f[x_i, x_0, x_1, \dots, x_{r-1}] \prod_{\substack{j=r \\ i \neq j}}^{n} \frac{x - x_j}{x_i - x_j} + \frac{f^{(n+1)}(\xi)}{(n+1)!} \prod_{i=r}^{n} (x - x_i) \tag{2.13}$$

**Proof.**

We know that

$$f[x, x_0, x_1, \dots, x_n] = \frac{f^{(n+1)}(\xi)}{(n+1)!}, \quad \text{Where } \xi \in \{x, x_0, x_1, \dots x_n\} \tag{2.14}$$

Using (2.14) in (2.1), we obtain equation (2.13)

**Corollary 2.1.2.** *Let* $x_0, x_1, \dots x_{r-1}$ *are 'r' numbers and* $x_r, x_{r+1}, \dots, x_n$ *are* $(n - r + 1)$ *distinct numbers in the interval* $[p, q]$, $r \leq n$, $x \in [p, q]$ *and* $f \in C^{n+1}[p, q]$, *then*

If $x \neq x_i, i = r, r+1, r+2, \dots, n$, *for some* $\xi \in \{x, x_0, x_1, \dots x_n\}$

(i).
$$f[x, x_0, x_1, \dots, x_{r-1}] = \prod_{j=r}^{n} (x - x_j) \sum_{i=r}^{n} \frac{f[x_i, x_0, x_1, \dots, x_{r-1}]}{x - x_i} \prod_{\substack{j=r \\ i \neq j}}^{n} \frac{1}{x_i - x_j} + \frac{f^{(n+1)}(\xi)}{(n+1)!} \prod_{i=r}^{n} (x - x_i) \tag{2.15}$$

(ii).
$$f[x, x_0, x_1, \dots, x_{r-1}] = \frac{\displaystyle\sum_{i=r}^{n} \frac{f[x_i, x_0, x_1, \dots, x_{r-1}]}{x - x_i} \prod_{\substack{j=r \\ i \neq j}}^{n} \frac{1}{x_i - x_j}}{\displaystyle\sum_{i=r}^{n} \frac{1}{x - x_i} \prod_{\substack{j=r \\ i \neq j}}^{n} \frac{1}{x_i - x_j}} + \frac{f^{(n+1)}(\xi)}{(n+1)!} \prod_{i=r}^{n} (x - x_i) \tag{2.16}$$

Generally, we construct Newton divided difference table to generate divided differences. But, we are unable to use this table to generate the divided differences in Theorem 2.1. So, we present new divided difference table to generate these divided differences.

### 2.1. New Divided Difference Table

In Newton divided difference table, divided differences of new entries in each column are determined by divided difference of two neighboring entries in the previous column. But, the procedure of new divided difference table is different from the Newton divided difference table. For example, consider the argument values $x_0, x_1$, $x_2 \ldots, x_6$ for the corresponding functional values $f_0, f_1, f_2, \ldots, f_6$. As a matter of convenience, we write $f_k = f(x_k)$. The procedure of New divided difference table is given in Table 1. The first order divided differences in the third column of the Table 1 are found by the sequence of evaluating $f[x_0, x_1]$, $f[x_0, x_2], \ldots$, The second order divided differences in the fourth column of the Table 1 are found by the sequence of evaluating $f[x_0, x_1, x_2]$, $f[x_0, x_1, x_3], \ldots$. Similarly, the sequences $f[x_0, x_1, x_2, x_3]$, $f[x_0, x_1, x_2, x_4], \ldots$ are evaluated for fifth column and the sequences $f[x_0, x_1, x_2, x_3, x_4]$, $f[x_0, x_1, x_2, x_3, x_5]$ are evaluated for sixth column and so on.

### 2.2. Combined form of Newton divided difference Table and New divided difference Table

Here, Newton divided difference table and new divided difference table are combined to produce a new combined form of divided difference table as shown in the Table 2. For example, consider the argument values $x_0, x_1, x_2 \ldots, x_6$ for the corresponding functional values $f_0, f_1, f_2, \ldots, f_6$. The table is divided into two parts and separated by a crossed line. The first part contains first four data follows the procedure of Newton divided difference table and remaining follows the procedure of new divided difference table.

The second part of the third column is divided difference of first and sixth data, first and seventh data

(i.e.). $f[x_0, x_5] = \dfrac{f_5 - f_0}{x_5 - x_0}$ and $f[x_0, x_6] = \dfrac{f_6 - f_0}{x_6 - x_0}$

The fourth column is divided difference of first and fifth and first and sixth

(i.e.). $f[x_0, x_1, x_5] = \dfrac{f[x_0, x_5] - f[x_0, x_1]}{x_5 - x_1}$ and $f[x_0, x_1, x_6] = \dfrac{f[x_0, x_6] - f[x_0, x_1]}{x_6 - x_1}$

The fifth column is divided difference of first and fourth and first and fifth and so on.

(i.e.). $f[x_0, x_1, x_2, x_5] = \dfrac{f[x_0, x_1, x_5] - f[x_0, x_1, x_2]}{x_5 - x_2}$ and $f[x_0, x_1, x_2, x_6] = \dfrac{f[x_0, x_1, x_6] - f[x_0, x_1, x_2]}{x_6 - x_2}$

Table 1. New divided difference table

| **x** | **y** | $\bar{\delta}^1$ | $\bar{\delta}^2$ | $\bar{\delta}^3$ | $\bar{\delta}^4$ |
|---|---|---|---|---|---|
| $x_0$ | $f_0$ | | | | |
| | | $f[x_0, x_1]$ | | | |
| $x_1$ | $f_1$ | | $f[x_0, x_1, x_2]$ | | |
| | | $f[x_0, x_2]$ | | $f[x_0, x_1, x_2, x_3]$ | |
| $x_2$ | $f_2$ | | $f[x_0, x_1, x_3]$ | | $f[x_0, x_1, x_2, x_3, x_4]$ |
| | | $f[x_0, x_3]$ | | $f[x_0, x_1, x_2, x_4]$ | |
| $x_3$ | $f_3$ | | $f[x_0, x_1, x_4]$ | | $f[x_0, x_1, x_2, x_3, x_5]$ |
| | | $f[x_0, x_4]$ | | $f[x_0, x_1, x_2, x_5]$ | |
| $x_4$ | $f_4$ | | $f[x_0, x_1, x_5]$ | | $f[x_0, x_1, x_2, x_3, x_6]$ |
| | | $f[x_0, x_5]$ | | $f[x_0, x_1, x_2, x_6]$ | |
| $x_5$ | $f_5$ | | $f[x_0, x_1, x_6]$ | | |
| | | $f[x_0, x_6]$ | | | |
| $x_6$ | $f_6$ | | | | |

Table 2. Combined form of Newton divided difference Table and New divided difference Table

| **x** | **y** | $\bar{\delta}^1$ | $\bar{\delta}^2$ | $\bar{\delta}^3$ | $\bar{\delta}^4$ |
|---|---|---|---|---|---|
| $x_0$ | $f_0$ | | | | |
| | | $f[x_0, x_1]$ | | | |
| $x_1$ | $f_1$ | | $f[x_0, x_1, x_2]$ | | |
| | | $f[x_1, x_2]$ | | $f[x_0, x_1, x_2, x_3]$ | |
| $x_2$ | $f_2$ | | $f[x_1, x_2, x_3]$ | | $f[x_0, x_1, x_2, x_3, x_4]$ |
| | | $f[x_2, x_3]$ | | $f[x_1, x_2, x_3, x_4]$ | |
| $x_3$ | $f_3$ | | $f[x_2, x_3, x_4]$ | | $f[x_0, x_1, x_2, x_3, x_5]$ |
| | | $f[x_3, x_4]$ | | $f[x_0, x_1, x_2, x_5]$ | |
| $x_4$ | $f_4$ | | $f[x_0, x_1, x_5]$ | | $f[x_0, x_1, x_2, x_3, x_6]$ |
| | | $f[x_0, x_5]$ | | $f[x_0, x_1, x_2, x_6]$ | |
| $x_5$ | $f_5$ | | $f[x_0, x_1, x_6]$ | | |
| | | $f[x_0, x_6]$ | | | |
| $x_6$ | $f_6$ | | | | |

### 2.3. Table of differences and divided differences by integer arguments.

Suppose the data are evenly spaced then, the construction of difference table is easier than new divided difference table. So, we have to change procedure of new combined form of divided difference table. Let $a, a+h, a+2h, \ldots, a+6h$ are the evenly spaced data, then their positional values from the top of the table are 0,1,2,…6. (For central difference, we can use $0, \pm 1, \pm 2, \ldots,$ as Shown in table 4). If we use the positional values, instead of using argument values, in combined form of divided difference table. Then, the divided differences in the first part of the Table 2 reduced to differences as shown in Table 3. But, from the 3rd column, first entry in each column of the table is represented by differences and factorials, this is because to find divided differences in the second part of the table and other entries contain only differences to reduce calculation burden. The divided differences in the second part of table are called as divided differences by integer arguments. For example $f_I[0,1,2,3,6]$, denote 4th order divided difference by integer arguments of 0,1,2,3,6 instead of $x_0, x_1, x_2, x_3, x_6$. The suffix $I$ denotes the divided difference is found by its integer arguments of the table.

Table 3. Table of differences and divided differences by Integer Arguments

| x | y | $\bar{\delta}^1$ | $\bar{\delta}^2$ | $\bar{\delta}^3$ | $\bar{\delta}^4$ |
|---|---|---|---|---|---|
| 0 | $f_0$ | | | | |
| | | $\frac{\Delta f_0}{1!}$ | | | |
| 1 | $f_1$ | | $\frac{\Delta^2 f_0}{2!}$ | | |
| | | $\Delta f_1$ | | $\frac{\Delta^3 f_0}{3!}$ | |
| 2 | $f_2$ | | $\Delta^2 f_1$ | | $\frac{\Delta^4 f_0}{4!}$ |
| | | $\Delta f_2$ | | $\Delta^3 f_1$ | |
| 3 | $f_3$ | | $\Delta^2 f_2$ | | $f_I[0,1,2,3,5]$ |
| | | $\Delta f_3$ | | $f_I[0,1,2,5]$ | |
| 4 | $f_4$ | | $f_I[0,1,5]$ | | $f_I[0,1,2,3,6]$ |
| | | $f_I[0,5]$ | | $f_I[0,1,2,6]$ | |
| 5 | $f_5$ | | $f_I[0,1,6]$ | | |
| | | $f_I[0,6]$ | | | |
| 6 | $f_6$ | | | | |

Table 4. Table of central differences and divided differences.

| x | y | $\bar{\delta}^1$ | $\bar{\delta}^2$ | $\bar{\delta}^3$ | $\bar{\delta}^4$ |
|---|---|---|---|---|---|
| −2 | $f_{-2}$ | | | | |
| | | $\delta f_{-2}$ | | | |
| -1 | $f_{-1}$ | | $\delta^2 f_{-2}$ | | |
| | | $\frac{\delta f_{-1}}{1!}$ | | $\frac{\delta^3 f_{-2}}{3!}$ | |
| 0 | $f_0$ | | $\frac{\delta^2 f_{-1}}{2!}$ | | $\frac{\delta^4 f_{-2}}{4!}$ |
| | | $\frac{\delta f_0}{1!}$ | | $\frac{\delta^3 f_{-1}}{3!}$ | |
| 1 | $f_1$ | | $\delta^2 f_0$ | | $f_I[0,-1,1,-2,3]$ |
| | | $\delta f_1$ | | $f_I[0,-1,1.3]$ | |
| 2 | $f_2$ | | $f_I[0,-1,3]$ | | $f_I[0,-1,1,-2,4]$ |
| | | $f_I[0,3]$ | | $f_I[0,-1,14]$ | |
| 3 | $f_3$ | | $f_I[0,-1,4]$ | | |
| | | $f_I[0,6]$ | | | |
| 4 | $f_4$ | | | | |

**Algorithm 2.1. (New Divided difference table).** Given the first two columns of the table, containing $x_0, x_1, \ldots, x_n$ and , correspondingly $f[x_0], f[x_1], \ldots, f[x_n]$, then the remaining entries are generated by the following steps,

Step 1: For $i$ = 1 to $r$ do
Step2 : For $j$ = 0 to $n$-$i$ do

Step3 : $$f[x_0, x_1, x_2, \ldots, x_{i-1}, x_{i+j}] = \frac{f[x_0, x_1, x_2, \ldots, x_{i-2}, x_{i+j}] - f[x_0, x_1, x_2, \ldots, x_{i-2}, x_{i-1}]}{x_{i+j} - x_{i-1}}$$

**Algorithm 2.2. (Combined form of divided difference table)** Given the first two columns of the table, containing $x_0, x_1, \ldots, x_n$ and , correspondingly $f[x_0], f[x_1], \ldots, f[x_n]$, then the remaining entries are generated by the following steps,

Step1 : For $i$ = 1 to $r$ do
Step2 : For $j$ = 0 to $n$-$i$ do

Step3 : If $j < r\text{-}i{+}1$ then $f\left[x_j, x_{j+1}, x_{j+2}, \ldots, x_{j+i}\right] = \dfrac{f\left[x_{j+1}, x_{j+2}, x_{j+3}, \ldots, x_{j+i}\right] - f\left[x_j, x_{j+1}, x_{j+2}, \ldots, x_{i+j-1}\right]}{x_{j+i} - x_j}$

Step4: Else $f\left[x_0, x_1, x_2, \ldots, x_{i-1}, x_{i+j}\right] = \dfrac{f\left[x_0, x_1, x_2, \ldots, x_{i-2}, x_{i+j}\right] - f\left[x_0, x_1, x_2, \ldots, x_{i-2}, x_{i-1}\right]}{x_{i+j} - x_{i-1}}$

**Algorithm 2.3. (Differences and divided differences by integer arguments)** Given the first two columns of the table, containing $0,1,2,\ldots,n$ and, correspondingly $f[x_0], f[x_1], \ldots, f[x_n]$, Let is a $d_{i,j}$ two dimensional array

Step1: For $i = 1$ to $r$ do
Step2: For $j = 0$ to $n\text{-}i$ do
Stpe3 : If $j < r\text{-}i{+}1$ then $d_{i,j} = d_{i-1,j+1} - d_{i-1,j}$

Step4: Else $d_{i,j} = \left( d_{i-1,j+1} - \dfrac{d_{i-1,0}}{i-1!} \right) \Big/ (j+1)$

**3. Interpolation formulas.**

**Theorem 3.1.** *Let $x_0, x_1, \ldots x_{r-1}$ are 'r' numbers and $x_r, x_{r+1}, \ldots, x_n$ are $(n-r+1)$ distinct numbers in the interval $[p,q]$, $r \le n$, $x \in [p,q]$ and $f \in C^{n+1}[p,q]$, then*

(i).
$$f(x) = \sum_{i=0}^{r-1} f[x_0, x_1, \ldots, x_i] \prod_{j=0}^{i-1} (x - x_j) + \prod_{i=0}^{r-1} (x - x_i) \sum_{i=r}^{n} f[x_i, x_0, x_1 \ldots, x_{r-1}] \prod_{\substack{j=r \\ i \ne j}}^{n} \frac{x - x_j}{x_i - x_j} + f\left[x, x_0, x_1, \ldots x_n\right] \prod_{i=0}^{n} (x - x_i) \tag{3.1}$$

(ii). *For $\xi \in \{x, x_0, x_1, \ldots x_n\}$*
$$f(x) = \sum_{i=0}^{r-1} f[x_0, x_1, \ldots, x_i] \prod_{j=0}^{i-1} (x - x_j) + \prod_{i=0}^{r-1} (x - x_i) \sum_{i=r}^{n} f[x_i, x_0, x_1, \ldots, x_{r-1}] \prod_{\substack{j=r \\ i \ne j}}^{n} \frac{x - x_j}{x_i - x_j} + \frac{f^{(n+1)}(\xi)}{(n+1)!} \prod_{i=0}^{n} (x - x_i) \tag{3.2}$$

**Proof.**

We know that

$$f[x, x_0, x_1, \ldots, x_{r-1}] = \frac{f[x, x_0, x_1, \ldots, x_{r-2}] - f[x_0, x_1, x_2 \ldots, x_{r-1}]}{x - x_{r-1}}$$

Rearranging this, we get

$$f[x, x_0, x_1, \ldots, x_{r-2}] = f[x_0, x_1, x_2 \ldots, x_{r-1}] + (x - x_{r-1}) f[x, x_0, x_1, \ldots, x_{r-1}] \tag{3.3}$$

Again, expanding $r-1^{th}$ divided difference of a polynomial as

$$f[x, x_0, x_1, \ldots, x_{r-3}] = f[x_0, x_1, x_2 \ldots, x_{r-2}] + (x - x_{r-2}) f[x_0, x_1, x_2 \ldots, x_{r-1}] + (x - x_{r-2})(x - x_{r-1}) f[x, x_0, x_1, \ldots, x_{r-1}] \tag{3.4}$$

Thus, repeating this, in general, for some '$m$', $m = 1,2,3,\ldots r$

$$\begin{aligned} f[x, x_0, x_1, \ldots, x_{r-m}] &= f[x_0, x_1, x_2 \ldots, x_{r-m+1}] + (x - x_{r-m+1}) f[x_0, x_1, x_2 \ldots, x_{r-m+2}] + \ldots \\ &+ (x - x_{r-m+1})(x - x_{r-m+2}) \ldots (x - x_{r-2}) f[x_0, x_1, x_2 \ldots, x_{r-1}] \\ &+ (x - x_{r-m+1})(x - x_{r-m+2}) \ldots (x - x_{r-1}) f[x, x_0, x_1, \ldots, x_{r-1}] \end{aligned} \quad (3.5)$$

When $m = r$ then,

$$\begin{aligned} f[x, x_0] &= f[x_0, x_1] + (x - x_1) f[x_0, x_1, x_2] + (x - x_1)(x - x_2) f[x_0, x_1, x_2, x_3] + \ldots \\ &+ (x - x_1)(x - x_2) \ldots (x - x_{r-2}) f[x_0, x_1, x_2 \ldots, x_{r-1}] \\ &+ (x - x_1)(x - x_2) \ldots (x - x_{r-1}) f[x, x_0, x_1, \ldots, x_{r-1}] \end{aligned}$$

After simplification, we get

$$\begin{aligned} f(x) &= f(x_0) + (x - x_0) f[x_0, x_1] + (x - x_0)(x - x_1) f[x_0, x_1, x_2] + \ldots \\ &+ (x - x_0)(x - x_1)(x - x_2) \ldots (x - x_{r-2}) f[x_0, x_1, x_2 \ldots, x_{r-1}] \\ &+ (x - x_0)(x - x_1) \ldots (x - x_{r-1}) f[x, x_0, x_1, \ldots, x_{r-1}] \end{aligned} \quad (3.6)$$

Using (2.1) in (3.6) and after simplification

$$\begin{aligned} &= f(x_0) + (x - x_0) f[x_0, x_1] + \ldots + (x - x_0)(x - x_1) \ldots (x - x_{r-2}) f[x_0, x_1, \ldots, x_{r-1}] \\ &+ \prod_{i=0}^{r-1} (x - x_i) \left( \sum_{i=r}^{n} f[x_i, x_0, x_1, \ldots, x_{r-1}] \prod_{\substack{i=r \\ i \neq j}}^{n} \frac{x - x_j}{x_i - x_j} \right) + f[x, x_0, x_1, \ldots x_n] \prod_{i=0}^{n} (x - x_i) \end{aligned} \quad (3.7)$$

After simplification, we obtain (3.1), we know that, there exists a number $\xi \in \{x, x_0, x_1, \ldots x_n\}$, then

$$f[x, x_0, x_1, \ldots, x_n] = \frac{f^{(n+1)}(\xi)}{(n+1)!} \quad (3.8)$$

Substituting (3.8) in equation (3.1), then we obtain Equation (3.2).

**Corollary 3.1.1 .** *Let* $x_0, x_1, \ldots x_{r-1}$ *are 'r' numbers and* $x_{r,} x_{r+1}, \ldots, x_n$ *are* $(n - r + 1)$ *distinct numbers in the interval* $[p, q]$, $r \leq n$ , $x \in [p, q]$ *and* $f \in C^{n+1}[p, q]$, *then*

Let $w_i^{(r)} = \prod_{\substack{j=r \\ i \neq j}}^{n} \frac{1}{x_i - x_j}$ , *For* $x \neq x_i$ , $i = r, r+1, r+2, \ldots, n$ *for some* $\xi \in \{x, x_0, x_1, \ldots x_n\}$

(i).

$$\begin{aligned} f(x) &= \sum_{i=1}^{r-1} f[x_0, x_1, \ldots, x_i] \prod_{j=0}^{i-1} (x - x_j) \\ &+ \prod_{i=0}^{n} (x - x_i) \sum_{i=r}^{n} \frac{f[x_i, x_0, x_1, \ldots, x_{r-1}]}{x - x_i} w_i^{(r)} + \prod_{i=0}^{n} (x - x_i) \frac{f^{(n+1)}(\xi)}{(n+1)!}) \end{aligned} \quad (3.9)$$

$$f(x)=\sum_{i=0}^{r-1} f[x_0,x_1,\ldots,x_i]\prod_{j=0}^{i-1}(x-x_j)$$

(ii). $$+\frac{\prod_{i=0}^{r-1}(x-x_i)\sum_{i=r}^{n}\frac{f[x_i,x_0,x_1,\ldots,x_{r-1}]}{x-x_i}w_i^{(r)}}{\sum_{i=r}^{n}\frac{1}{x-x_i}w_i^{(r)}}+\frac{f^{(n+1)}(\xi)}{(n+1)!}\prod_{i=0}^{n}(x-x_i) \qquad (3.10)$$

If we put $r=n$ and $r=0$ in Theorem 3.1, we obtain Newton's general interpolation formula and Lagrange's interpolation formula respectively. If we put $r=0$ in (3.10), we obtain bary centric interpolation formula.

**Corollary 3.1.2.** *Let* $x_0,x_1,x_2\ldots x_n$ *are* $(n+1)$ *distinct numbers in the interval* $[p,q]$, *spaced equally i.e* $x_i=x_0+ih,(i=0,1,2,\ldots n)$, $h\neq 0$ *and* $f\in C^{n+1}[p,q]$, *if* $x=x_0+sh,$ *then for some* $\xi\in\{x,x_0,x_1,\ldots x_n\}$

$$f(x_0+sh)=f_0+\sum_{i=1}^{r-1}\Delta^i f_0\frac{s^{(i)}}{i!}+s^{(r)}\sum_{i=r}^{n}f_I[i,0,1,\ldots,r-1]\prod_{\substack{j=r\\ i\neq j}}^{n}\frac{s-j}{i-j}+\frac{h^{n+1}s^{(n+1)}}{(n+1)!}f^{(n+1)}(\xi) \qquad (3.11)$$

Where $\Delta^i f_0$ is the $i^{th}$ difference and $f_I[i,0,1,\ldots,r-1]$ is the $r^{th}$ divided difference by integer arguments

**Corollary 3.1.3.** *Let* $x_0,x_1,x_2\ldots x_n$ *are* $(n+1)$ *distinct numbers in the interval* $[p,q]$, $x\in[p,q]$ *spaced equally i.e,* $x_i=x_0-ih,(i=0,1,2,\ldots n)$, $h\neq 0$ *and* $f\in C^{n+1}[p,q]$, *if* $x=x_0+sh,$ *then for some* $\xi\in\{x,x_0,x_1,\ldots x_n\}$

$$f(x_0+sh)=f_0+\sum_{i=1}^{r-1}\frac{\nabla^i f_0 s^{(-i)}}{i!}+s^{(-r)}(-1)^{n-r}\sum_{i=r}^{n}f_I[-i,0,-1,\ldots,-(r-1)]\prod_{\substack{j=r\\ i\neq j}}^{n}\frac{s+j}{i-j}+\frac{h^{n+1}s^{(-n-1)}f^{(n+1)}(\xi)}{(n+1)!} \qquad (3.12)$$

Where $\nabla^i f_0$ is $i^{th}$ backward difference and $f_I[-i,0,-1,\ldots,\overline{-r-1}]$ is $r^{th}$ divided difference by integer arguments

**Corollary 3.1.4.** *Let* $x_{-m}\ldots x_{-2},x_{-1},x_0,x_1,x_2\ldots x_n$ *are the* $(n+m+1)$ *distinct numbers in the interval* $[p,q]$, $x\in[p,q]$ *spaced equally, i.e.* $x_i=x_0+ih,(i=-m,\ldots,-2,-1,0,1,2,\ldots n)$, $h\neq 0$ *and* $f\in C^{n+m+1}[p,q]$, *if* $x=x_0+sh,$ *then for some* $\xi\in\{x,x_{-m}\ldots x_{-2},x_{-1,}x,x_0,x_1,\ldots x_n\}$

(i). $$f(x_0+sh)=f_0+\frac{\delta^1 f_0}{1!}s^{(1)}+\frac{\delta^2 f_{-1}}{2!}s^{(2)}+\frac{\delta^3 f_{-1}}{3!}(s+1)^{(3)}+\ldots+\frac{\delta^{2r}f_{-r}}{2r!}(s+r-1)^{(2r)}$$
$$+(s+r)^{(2r+1)}\sum_{i=\pm(r+1)}^{-m,+n}f_I[i,0,1,-1,2,-2,\ldots r,-r]\prod_{\substack{j=\pm(r+1)\\ i\neq j}}^{-m,+n}\frac{s-j}{i-j}+\frac{h^{n+m+1}f^{(n+m+1)}(\xi)}{(n+m+1)!}\prod_{i=-m}^{n}(s-i) \qquad (3.13)$$

(ii). $$f(x_0+sh)=f_0+\frac{\delta^1 f_{-1}}{1!}s+\frac{\delta^2 f_{-1}}{2!}(s+1)^{(2)}+\frac{\delta^3 f_{-2}}{3!}(s+1)^{(3)}+\ldots+\frac{\delta^{2r}f_{-r}}{2r!}(s+r)^{(2r)}$$
$$+(s+r)^{(2r+1)}\sum_{i=\pm(r+1)}^{-m,n}f_I[i,0,-1,1,-2,2,\ldots,-r,r]\prod_{\substack{j=\pm(r+1)\\ i\neq j}}^{-m,n}\frac{s-j}{i-j}+\frac{h^{n+m+1}f^{(n+m+1)}(\xi)}{(n+m+1)!}\prod_{i=-m}^{n}(s-i) \qquad (3.14)$$

(3.13) and (3.14) are new central difference forward and backward difference formula respectively.

**Various forms of Central difference formulas**

Let $\Theta(s) = (s+r)^{(2r+1)} \sum_{i=\pm(r+1)}^{-m,+n} f_I[i,0,-1,1,-2,2,\ldots,-r,r] \prod_{\substack{j=\pm(r+1)\\ i\neq j}}^{-m,+n} \frac{s-j}{i-j}$, $R_{m+n+1} = \frac{h^{m+n+1} f^{(m+n+1)}(\xi)}{(m+n+1)!} \prod_{i=-m}^{n}(s-i)$

and $\xi \in \{x, x_{-m} \ldots x_{-2}, x_{-1,} x, x_0, x_1, \ldots x_n\}$

Using Corollary (3.13) and (3.14), we can find other interpolation formulas using central differences, few of them listed below

(i).
$$f(x_0+sh) = f_0 + \frac{\delta^1 f(0)+\delta^1 f(-1)}{2}\frac{s}{1!} + \frac{s^2}{2!}\delta^2 f(-1) + \ldots + \frac{\delta^{2r-1}f(-r+1)+\delta^{2r-1}f(-r)}{2}\frac{(s+r-1)^{(2r-1)}}{2r-1!}$$
$$+\delta^{2r}f(-r)\frac{s(s+r-1)^{(2r-1)}}{2r!} + \Theta(s) + R_{n+m+1} \quad (3.15)$$

(ii).
$$f(x_0+sh) = \frac{f_0+f_1}{2} + \frac{(s-\frac{1}{2})}{1!}\delta f(0) + \frac{s(s-1)}{2!}\frac{\delta^2 f(-1)+\delta^2 f(0)}{2} + \frac{(s-\frac{1}{2})s(s-1)}{3!}\delta^3 f(-1)$$
$$+\frac{(s+1)s(s-1)(s-2)}{4!}\frac{\delta^4 f(-2)+\delta^4 f(-1)}{2} + \ldots + \frac{(s+r-1)^{(2r)}}{2r!}\frac{\delta^{2r}f(-r)+\delta^{2r}f(-r+1)}{2}$$
$$-\frac{(s+r-1)^{(2r)}}{2r!}\frac{\delta^{2r+1}f(-r)}{2} + \Theta(s) + R_{n+m+1} \quad (3.16)$$

(iii).
$$f(x_0+sh) = f_0 t + \frac{(t+1)^{(3)}}{3!}\delta^2 f(-1) + \ldots + \frac{\delta^{2r-2}f(-r+1)}{2r-1!}(t+r-1)^{(2r-1)} + f_1 s + \frac{(s+1)^{(3)}}{3!}\delta^2 f(0)$$
$$+\ldots+\frac{\delta^{2r-2}f(-r+2)}{2r-1!}(s+r-1)^{(2r-1)} + \frac{\delta^{2r}f(-r)}{2r!}(s+r-1)^{(2r)} + \Theta(s) + R_{n+m+1} \quad (3.17)$$

(iv).
$$f(x_0+sh) = f_0 + \delta f(0)\frac{(s+1)^{(2)}}{2!} - \frac{\delta f(-1)}{2!}s^{(2)} + \ldots + \delta^{2r-1}f(-r+1)\frac{(s+r)^{(2r)}}{2r!}$$
$$-\frac{\delta^{2r-1}f(-r)}{2r!}(s+r-1)^{(2r)} + \Theta(s) + R_{n+m+1} \quad (3.18)$$

Where $t = 1 - s$, Equation (3.15), (3.16), (3.17) and (3.18) are new modified form of Striling's, Bessels's, Everett's and Steffensen's central difference formulas respectively.

**Approximation by New Interpolation formula**

The general structure of the new interpolation formula can be written in the following form

$$f(x) = N(x) + (x-x_0)(x-x_1)(x-x_2)\ldots(x-x_{r-1})\lambda(x) + E(x) \quad (3.19)$$

Where $f_0, f_1, f_2, \ldots, f_n$ are the functional values of the function *f*, for the distinct arguments $x_0, x_1, x_2 \ldots x_n$ and $f \in C^{n+1}[p,q]$. *N*(*x*) is Newton's interpolation formula up to certain order and $\lambda(x)$ is the divided difference polynomial. To approximate new interpolation formula, we can replace $\lambda(x)$ by a suitable approximation $\Pi(x)$. We can use least squares or any other methods to replace $\lambda(x)$ by $\Pi(x)$.

i.e. $f(x) = N(x) + (x-x_0)(x-x_1)(x-x_2)\ldots(x-x_{r-1})\Pi(x)$ (3.20)

$\Pi(x)$ is an approximation of $r^{th}$ divided difference polynomial, obtained from the new divided difference table. The replacement of $\Pi(x)$ by the suitable approximation is shown through an example. Here, we use Method of least squares to replace $\Pi(x)$. The actual function taken for comparison is $f(x) = e^x(1+x) + x\sin x$ and the argument values $x$ and corresponding functional values $y$ are given in the following table.

Table 5

| $x$ | 1.0 | 1.25 | 1.5 | 1.75 | 2.0 | 2.25 | 2.5 | 2.75 | 3.00 |
|---|---|---|---|---|---|---|---|---|---|
| $y$ | 6.2780346 | 9.0395024 | 12.7004652 | 17.5471328 | 23.9857632 | 32.5858062 | 44.1349092 | 59.7094373 | 80.7655077 |

In general, former formulas of interpolation based on differences are compared by their speed of convergence. So, we consider the 4th the degree polynomial, generated by Newton forward and backward, Gauss forward and backward, Stirling's, Bessel's, Everett's and Steffensen's interpolation formulas in equations (3.21), (3.22), (3.23), (3.24), (3.25), (3.26), (3.27) and (3.28) respectively. The errors of various interpolation formulas with actual functional values on (0.85, 3.15) are shown in the Table 6.

$$f(1+sh) = 6.278035 + 2.761468s + 0.449747s(s-1) + 0.047702s(s-1)(s-2) + 0.005002s(s-1)(s-2)(s-3) \quad (3.21)$$

Where $s = (x-1)/0.25$

$$f(3+sh) = 80.765508 + 21.056070s + 2.740771s(s+1) + 0.242686s(s+1)(s+2) + 0.015823s(s+1)(s+2)(s+3) \quad (3.22)$$

Where $s = (x-3)/0.25$

$$f(2+sh) = 23.985763 + 8.600043s + 1.080706s(s-1) + 0.131275s(s^2-1) + 0.009092s(s^2-1)(s-2) \quad (3.23)$$

$$f(2+sh) = 23.985763 + 6.438630s + 1.080706s(s+1) + 0.094908s(s^2-1) + 0.009092s(s^2-1)(s+2) \quad (3.24)$$

$$f(2+sh) = 23.985763 + 7.519337s + 1.080706s^2 + 0.113091s(s^2-1) + 0.009092s^2(s^2-1) \quad (3.25)$$

$$f(2+sh) = 28.285785 + 8.600043(s-0.5) + 1.277618s(s-1) + 0.131275s(s-0.5)(s-1) + 0.010561s(s^2-1)(s-2) \quad (3.26)$$

$$f(2+sh) = 23.985763t + 0.360235t(t^2-1) + 32.585806s + 0.491510s(s^2-1) + 0.009092s(s^2-1)(s-2) \quad (3.27)$$

$$f(2+sh) = 23.985763 + 4.300022(s+1)s - 3.219315s(s-1) + 0.032819s(s^2-1)(s+2) - 0.023727s(s^2-1)(s-2) \quad (3.28)$$

Where $t = 1-s$ and $s = (x-2)/0.25$

Similarly, The new modified formulas (of Newton forward, Newton backward, Gauss forward, Gauss backward, Stirling's, Bessel's, Everett's and Steffensen's interpolation formulas) are used to calculate $4^{th}$ degree polynomial and they are given in equations (3.29), (3.30), (3.31), (3.32), (3.33), (3.34), (3.35) and (3.36) respectively. The divided difference polynomial in (3.19) is replaced by the function $\theta(s)$, where $\theta(s)$ is a linear function found by using method of least squares, on $3^{rd}$ order divided differences generated in the new divided difference table. The errors of various new interpolation formulas and their actual functional values on (0.85, 3.15) are shown in the Table 7.

$$f(1+sh) = 6.278035 + 2.761468s + 0.449747s(s-1) + s(s-1)(s-2)\theta(s) \quad (3.29)$$

Where $\theta(s) = 0.006633s + 0.026390$

$$f(3+sh) = 80.765508 + 21.056070s + 2.740771s(s+1) + s(s+1)(s+2)\theta(s) \quad (3.30)$$

Where $\theta(s) = 0.013078s + 0.279710$

$$f(2+sh) = 23.985763 + 8.600043s + 1.080706s(s-1) + s(s^2-1)\theta(s) \quad (3.31)$$

$$f(2+sh) = 23.985763 + 6.438630s + 1.080706s(s+1) + s(s^2-1)\theta(s) \quad (3.32)$$

$$f(2+sh) = 23.985763 + 7.519337s + 1.080706s^2 + s(s^2-1)\theta(s) \quad (3.33)$$

$$f(2+sh) = 28.285785 + 8.600043(s-0.5) + 1.277618s(s-1) - 0.196912s(s-1) + s(s^2-1)\theta(s) \quad (3.34)$$

$$f(2+sh) = 23.985763t + +32.585806s + 1.080706s(s-1) + s(s^2-1)\theta(s) \quad (3.35)$$

$$f(2+sh) = 23.985763 + 4.300022(s+1)s - 3.219315s(s-1) + s(s^2-1)\theta(s) \quad (3.36)$$

Where $\theta(s) = 0.009269s + 0.116079$

Table 6. Errors of Equations (3.21) to (3.28) with actual function values on various values of $x$

| x | Newton Forward Formula | Newton Backward Formula | Central difference: Forward Formula | Backward Formula | Strirling's formula | Bessel's formula | Everett's formula | Steffensen's formula |
|---|---|---|---|---|---|---|---|---|
| 0.85 | 1.19e-02 | 7.39e+00 | 6.88e-01 | 6.88e-01 | 6.88e-01 | 1.59e+00 | 6.88e-01 | 6.88e-01 |
| 0.90 | 5.81e-03 | 6.33e+00 | 5.41e-01 | 5.41e-01 | 5.41e-01 | 1.30e+00 | 5.41e-01 | 5.41e-01 |
| 1.00 | 0.00e+00 | 4.56e+00 | 3.18e-01 | 3.18e-01 | 3.18e-01 | 8.47e-01 | 3.18e-01 | 3.18e-01 |
| 1.10 | -1.07e-03 | 3.20e+00 | 1.73e-01 | 1.73e-01 | 1.73e-01 | 5.27e-01 | 1.73e-01 | 1.73e-01 |
| 1.15 | -8.23e-04 | 2.65e+00 | 1.23e-01 | 1.23e-01 | 1.23e-01 | 4.07e-01 | 1.23e-01 | 1.23e-01 |
| 1.25 | 0.00e+00 | 1.77e+00 | 5.50e-02 | 5.50e-02 | 5.50e-02 | 2.31e-01 | 5.50e-02 | 5.50e-02 |
| 1.35 | 4.33e-04 | 1.14e+00 | 1.92e-02 | 1.92e-02 | 1.92e-02 | 1.20e-01 | 1.92e-02 | 1.92e-02 |
| 1.50 | 0.00e+00 | 5.26e-01 | 0.00e+00 | 0.00e+00 | 0.00e+00 | 3.53e-02 | 3.47e-18 | 0.00e+00 |
| 1.65 | -4.54e-04 | 2.03e-01 | -1.34e-03 | -1.34e-03 | -1.34e-03 | 5.38e-03 | -1.34e-03 | -1.34e-03 |
| 1.75 | 0.00e+00 | 9.10e-02 | 0.00e+00 | 0.00e+00 | 0.00e+00 | 0.00e+00 | 0.00e+00 | 0.00e+00 |
| 1.85 | 9.21e-04 | 3.19e-02 | 7.03e-04 | 7.03e-04 | 7.03e-04 | -7.64e-04 | 7.03e-04 | 7.03e-04 |
| 2.00 | 1.73e-18 | 0.00e+00 | 0.00e+00 | 0.00e+00 | 0.00e+00 | 0.00e+00 | 0.00e+00 | 0.00e+00 |
| 2.10 | -7.04e-03 | -2.88e-03 | -6.76e-04 | -6.76e-04 | -6.76e-04 | 1.14e-04 | -6.76e-04 | -6.76e-04 |
| 2.15 | -1.46e-02 | -2.22e-03 | -7.39e-04 | -7.39e-04 | -7.39e-04 | 5.11e-05 | -7.39e-04 | -7.39e-04 |
| 2.25 | -4.31e-02 | 0.00e+00 | 0.00e+00 | 0.00e+00 | 0.00e+00 | 0.00e+00 | 0.00e+00 | 0.00e+00 |
| 2.35 | -9.93e-02 | 1.17e-03 | 1.50e-03 | 1.50e-03 | 1.50e-03 | 3.15e-04 | 1.50e-03 | 1.50e-03 |
| 2.50 | -2.71e-01 | 0.00e+00 | 0.00e+00 | 0.00e+00 | 0.00e+00 | 0.00e+00 | 0.00e+00 | 0.00e+00 |
| 2.65 | -6.15e-01 | -1.23e-03 | -2.38e-02 | -2.38e-02 | -2.38e-02 | -1.06e-02 | -2.38e-02 | -2.38e-02 |
| 2.75 | -9.93e-01 | 0.00e+00 | -7.05e-02 | -7.05e-02 | -7.05e-02 | -3.53e-02 | -7.05e-02 | -7.05e-02 |
| 2.85 | -1.54e+00 | 2.51e-03 | -1.63e-01 | -1.63e-01 | -1.63e-01 | -8.87e-02 | -1.63e-01 | -1.63e-01 |
| 3.00 | -2.78e+00 | 0.00e+00 | -4.44e-01 | -4.44e-01 | -4.44e-01 | -2.67e-01 | -4.44e-01 | -4.44e-01 |
| 3.10 | -3.99e+00 | -1.93e-02 | -7.80e-01 | -7.80e-01 | -7.80e-01 | -4.95e-01 | -7.80e-01 | -7.80e-01 |
| 3.15 | -4.74e+00 | -4.00e-02 | -1.01e+00 | -1.01e+00 | -1.01e+00 | -6.56e-01 | -1.01e+00 | -1.01e+00 |

Table 7. Errors of Equations (3.29) to (3.36) with actual functional values on various values of $x$

| x | Modified Newton Forward Formula | Modified Newton Backward Formula | Cenntral difference: Modified Forward Formula | Modified Backward Formula | Modified Strirling's Formula | Modified Bessel's Formula | Modified Everett's Formula | Modified Steffensen's Formula |
|---|---|---|---|---|---|---|---|---|
| 0.85 | 3.00e-02 | 1.71e+00 | 4.87e-01 | 4.87e-01 | 4.87e-01 | 4.87e-01 | 4.87e-01 | 4.87e-01 |
| 0.90 | 1.52e-02 | 1.31e+00 | 3.62e-01 | 3.62e-01 | 3.62e-01 | 3.62e-01 | 3.62e-01 | 3.62e-01 |
| 1.00 | 0.00e+00 | 6.93e-01 | 1.82e-01 | 1.82e-01 | 1.82e-01 | 1.82e-01 | 1.82e-01 | 1.82e-01 |
| 1.10 | -3.24e-03 | 2.77e-01 | 7.19e-02 | 7.19e-02 | 7.19e-02 | 7.19e-02 | 7.19e-02 | 7.19e-02 |
| 1.15 | -2.61e-03 | 1.29e-01 | 3.70e-02 | 3.70e-02 | 3.70e-02 | 3.70e-02 | 3.70e-02 | 3.70e-02 |
| 1.25 | 0.00e+00 | -6.88e-02 | -3.91e-03 | -3.91e-03 | -3.91e-03 | -3.91e-03 | -3.91e-03 | -3.91e-03 |
| 1.35 | 1.78e-03 | -1.68e-01 | -1.86e-02 | -1.86e-02 | -1.86e-02 | -1.86e-02 | -1.86e-02 | -1.86e-02 |
| 1.50 | 0.00e+00 | -1.97e-01 | -1.58e-02 | -1.58e-02 | -1.58e-02 | -1.58e-02 | -1.58e-02 | -1.58e-02 |
| 1.65 | -5.61e-03 | -1.51e-01 | -5.02e-03 | -5.02e-03 | -5.02e-03 | -5.02e-03 | -5.02e-03 | -5.02e-03 |
| 1.75 | -8.49e-03 | -1.06e-01 | 0.00e+00 | 0.00e+00 | 0.00e+00 | 0.00e+00 | 0.00e+00 | 0.00e+00 |
| 1.85 | -7.79e-03 | -6.21e-02 | 1.81e-03 | 1.81e-03 | 1.81e-03 | 1.81e-03 | 1.81e-03 | 1.81e-03 |
| 2.00 | 5.19e-03 | -1.28e-02 | 0.00e+00 | 0.00e+00 | 0.00e+00 | 0.00e+00 | 0.00e+00 | 0.00e+00 |
| 2.10 | 2.41e-02 | 5.54e-03 | -1.70e-03 | -1.70e-03 | -1.70e-03 | -1.70e-03 | -1.70e-03 | -1.70e-03 |
| 2.15 | 3.68e-02 | 1.05e-02 | -1.93e-03 | -1.93e-03 | -1.93e-03 | -1.93e-03 | -1.93e-03 | -1.93e-03 |
| 2.25 | 6.77e-02 | 1.33e-02 | 0.00e+00 | 0.00e+00 | 0.00e+00 | 0.00e+00 | 0.00e+00 | 0.00e+00 |
| 2.35 | 1.03e-01 | 9.43e-03 | 5.85e-03 | 5.85e-03 | 5.85e-03 | 5.85e-03 | 5.85e-03 | 5.85e-03 |
| 2.50 | 1.47e-01 | 0.00e+00 | 2.01e-02 | 2.01e-02 | 2.01e-02 | 2.01e-02 | 2.01e-02 | 2.01e-02 |
| 2.65 | 1.42e-01 | -3.45e-03 | 2.78e-02 | 2.78e-02 | 2.78e-02 | 2.78e-02 | 2.78e-02 | 2.78e-02 |
| 2.75 | 7.99e-02 | 0.00e+00 | 1.40e-02 | 1.40e-02 | 1.40e-02 | 1.40e-02 | 1.40e-02 | 1.40e-02 |
| 2.85 | -6.32e-02 | 5.46e-03 | -3.36e-02 | -3.36e-02 | -3.36e-02 | -3.36e-02 | -3.36e-02 | -3.36e-02 |
| 3.00 | -5.15e-01 | 0.00e+00 | -2.22e-01 | -2.22e-01 | -2.22e-01 | -2.22e-01 | -2.22e-01 | -2.22e-01 |
| 3.10 | -1.05e+00 | -3.48e-02 | -4.75e-01 | -4.75e-01 | -4.75e-01 | -4.75e-01 | -4.75e-01 | -4.75e-01 |
| 3.15 | -1.41e+00 | -7.02e-02 | -6.57e-01 | -6.57e-01 | -6.57e-01 | -6.57e-01 | -6.57e-01 | -6.57e-01 |

Newton forward and backward formulas give fair accuracy near the beginning and end of the Table 6 respectively. But, at the central zone of Table 6, central difference formulas give fair accuracy than Newton forward and backward formulas. But, comparing the results of corresponding former formulas in Table 6 to new modified formulas in Table 7, it is clear that the new modified formulas give fair accuracy than the former formulas. Through out the interval (0.85, 3.15) new modified formulas possess good accuracy. However, comparing the results of new modified formula in Table 7, forward formula, backward formula and central difference formulas gives better accuracy in their respective(forward, backward and central) positions.

Thus, the new interpolation formulas posses some extra advantages then the former interpolation formulas, we can restrict the divided differences up to certain order, instead of finding all divided differences in the Newton divided difference formula. Also, we can replace a part of the formula by a suitable approximation. The table 8 shows the comparisons of number of various operations between Newton‘s formula, Lagrange’s formula and new interpolation formula. Newton’s formula needs more number of divisions but requires less number of multiplications than Lagrange formula whereas Lagrange formula needs more number of multiplications but requires less number of divisions than Newton’s formula. But, in the new formula, using values of $r$ we can choose which operation is to be used maximum. Both Newton’s and Lagrange’s interpolation formulas are the special case of new interpolation formula.

Table 8. Comparisons of Total number of different operations of various interpolation formulas.

| Operations | Newton Divided difference formula | Lagrange Interpolation formula | New Interpolation formula, Theorem 3.1 $0<r<n$ |
|---|---|---|---|
| Additions | $n$ | $n$ | $n$ |
| Subtractions | $3n(n+1)/2$ | $2n(n+1)$ | $n(n+1)+(n-r)(n-r+1)+r(r+1)/2$ |
| Multiplications | $n(n+1)/2$ | $(2n-1)(n+1)$ | $(2(n-r)-1)(n-r+1)+r(r+1)/2$ |
| Division | $n(n+1)/2$ | $n+1$ | $n(n+1)/2-(n-r)(n-r+1)/2+n-r+1$ |

## 4. Formulas of Numerical Differentiation

In this section, we present two new different types of formulas to evaluate higher order derivatives in the following form.

1. Recursive formula approach
2. Linear combination of divided differences

### 4.1. Recursive formula for $k^{th}$ order differentiation

**Theorem 4.1.** *Let* $x, x_0, x_1, \dots x_n$ *are* $(n+1)$ *distinct numbers in the interval* $[p,q]$ , $k \in W$ *and* $f \in C^{n+k+1}[p,q]$ *then*

(i).
$$\frac{f^{(k)}(x)}{k!}\rho_0(x)+\frac{f^{(k-1)}(x)}{k-1!}\rho_1(x)+\dots+\frac{f(x)}{0!}\rho_k(x) = \sum_{i=0}^{n}\frac{f(x_i)}{(x_i-x)^k}\prod_{\substack{j=0\\i\neq j}}^{n}\frac{x-x_j}{x_i-x_j}+f[\underbrace{x,\dots x}_{k+1\ times},x_0,x_1,x_2\dots x_n]\prod_{i=0}^{n}(x-x_i) \tag{4.1}$$

(ii).
$$\frac{f^{(k)}(x)}{k!}\rho_0(x)+\frac{f^{(k-1)}(x)}{k-1!}\rho_1(x)+\dots+\frac{f(x)}{0!}\rho_k(x) = \sum_{i=0}^{n}\frac{f(x_i)}{(x_i-x)^k}\prod_{\substack{j=0\\i\neq j}}^{n}\frac{x-x_j}{x_i-x_j}+\frac{f^{(k+n+1)}(\xi)}{(k+n+1)!}\prod_{i=0}^{n}(x-x_i) \tag{4.2}$$

Where $\rho_0(x)=1, \rho_m(x)=\sum_{i=0}^{n}\frac{1}{(x_i-x)^m}\prod_{\substack{j=0\\i\neq j}}^{n}\frac{x-x_j}{x_i-x_j}$ , $m=1,2,3,\dots k$ and $\xi \in \{x, x_0, x_1, \dots x_n\}$

**Proof.**

Applying Theorem (2.1) for $f[\underbrace{x,\ldots,x}_{k+1\,times}]$

$$f[\underbrace{x,\ldots,x}_{k+1\,times}] = \sum_{i=0}^{n} f[\underbrace{x,\ldots,x}_{k\,times},x_i]\prod_{\substack{j=0\\i\neq j}}^{n}\frac{x-x_j}{x_i-x_j} + f[\underbrace{x,\ldots x}_{k+1\,times},x_0,x_1,x_2\ldots x_n]\prod_{i=0}^{n}(x-x_i) \tag{4.3}$$

$$\frac{f^{(k)}(x)}{k!} = \sum_{i=0}^{n} f[\underbrace{x,\ldots,x}_{k\,times},x_i]\prod_{\substack{j=0\\i\neq j}}^{n}\frac{x-x_j}{x_i-x_j} + f[\underbrace{x,\ldots x}_{k+1\,times},x_0,x_1,x_2\ldots x_n]\prod_{i=0}^{n}(x-x_i) \tag{4.4}$$

We know that

$$f[x_i,\underbrace{x,\ldots,x}_{k\,times}] = \frac{f[x_i,\underbrace{x,\ldots,x}_{k-1\,times}] - f[\underbrace{x,\ldots,x}_{k\,times}]}{x_i - x}$$

$$= \frac{f[x_i,\underbrace{x,\ldots,x}_{k-1\,times}]}{x_i - x} - \frac{f^{(k-1)}(x)}{k-1!(x_i-x)}$$

$$= -\frac{f^{(k-1)}(x)}{k-1!(x_i-x)} - \frac{f^{(k-2)}(x)}{k-2!(x_i-x)^2} + \frac{f[x_i,\underbrace{x,\ldots,x}_{k-2\,times}]}{(x_i-x)^2}$$

Preceding this, we get

$$= -\frac{f^{(k-1)}(x)}{k-1!(x_i-x)} - \frac{f^{(k-2)}(x)}{k-2!(x_i-x)^2} - \ldots - \frac{f(x)}{0!(x_i-x)^k} + \frac{f(x_i)}{0!(x_i-x)^k} \tag{4.5}$$

Substituting (4.5) in (4.4) we get

$$\begin{aligned}\frac{f^{(k)}(x)}{k!} &= \sum_{i=0}^{n}\left\{-\frac{f^{(k-1)}(x)}{k-1!(x_i-x)} - \frac{f^{(k-2)}(x)}{k-2!(x_i-x)^2} - \ldots - \frac{f(x)}{0!(x_i-x)^k} + \frac{f(x_i)}{0!(x_i-x)^k}\right\}\prod_{\substack{j=0\\i\neq j}}^{n}\frac{x-x_j}{x_i-x_j}\\ &\quad + f[\underbrace{x,\ldots x}_{k+1\,times},x_0,x_1,x_2\ldots x_n]\prod_{i=0}^{n}(x-x_i)\end{aligned} \tag{4.6}$$

After simplification, we get

$$\begin{aligned}&\frac{f^{(k)}(x)}{k!} + \frac{f^{(k-1)}(x)}{k-1!}\sum_{i=0}^{n}\frac{1}{(x_i-x)}\prod_{\substack{j=0\\i\neq j}}^{n}\frac{x-x_j}{x_i-x_j} + \ldots + \frac{f(x)}{0!}\sum_{i=0}^{n}\frac{1}{(x_i-x)^k}\prod_{\substack{j=0\\i\neq j}}^{n}\frac{x-x_j}{x_i-x_j}\\ &= \sum_{i=0}^{n}\frac{f(x_i)}{(x_i-x)^k}\prod_{\substack{j=0\\i\neq j}}^{n}\frac{x-x_j}{x_i-x_j} + f[\underbrace{x,\ldots x}_{k+1\,times},x_0,x_1,x_2\ldots x_n]\prod_{i=0}^{n}(x-x_i)\end{aligned} \tag{4.7}$$

Put $\rho_m(x)=\sum_{i=0}^{n}\frac{1}{(x_i-x)^m}\prod_{\substack{j=0\\ i\neq j}}^{n}\frac{x-x_j}{x_i-x_j}$ , $m=1,2,\ldots,k$ and $\rho_0(x)=1$ in (4.7) (4.8)

Thus, we obtain equation (4.1). Also, We know that

$$f[\underbrace{x,\ldots x}_{k+1\ times},x_0,x_1,x_2\ldots x_n]=\frac{f^{(k+n+1)}(\xi)}{(k+n+1)!}\text{, for some } \xi\in\{x,x_0,x_1,\ldots x_n\} \tag{4.9}$$

Substituting (4.9) in (4.1) we obtain (4.2).

**Corollary 4.1.1.** *Let* $x,x_0,x_1,\ldots x_n$ *are distinct numbers in the interval* $[p,q]$, $k\in W$ *and* $f\in C^{n+k+1}[p,q]$, *then*

$$\begin{aligned}&\frac{f^{(k)}(x)}{k!}\prod_{i=0}^{n}(x_i-x)^k\phi_0(x)+\prod_{i=0}^{n}(x_i-x)^{k-1}\frac{f^{(k-1)}(x)}{k-1!}\phi_1(x)+\ldots+\frac{f(x)}{0!}\phi_k(x)\\&=\sum_{i=0}^{n}f(x_i)\prod_{\substack{j=0\\ i\neq j}}^{n}\frac{(x_j-x)^{k+1}}{x_i-x_j}+\frac{f^{(k+n+1)}(\xi)}{(k+n+1)!}\prod_{i=0}^{n}(x_i-x)^{k+1}\end{aligned} \tag{4.10}$$

Where $\phi_m(x)=\sum_{i=0}^{n}\prod_{\substack{j=0\\ i\neq j}}^{n}\frac{(x_j-x)^{m+1}}{x_i-x_j}$, $m=0,1,2,3,\ldots k$ and for some $\xi\in\{x,x_0,x_1,\ldots x_n\}$

The following three notes are very handy to derive various forms of Equation 4.2 when the data are equally spaced.

**Note 4.1.** Suppose the data are spaced equally as $a,a+h,a+2h,\ldots,a+nh$ , then value of the $\rho_m(a)$ is calculated as follow as

$$\begin{aligned}\prod_{\substack{j=1\\ i\neq j}}^{n}\frac{j}{j-i}&=\frac{1\times2\times\ldots\times(i-1)\times(i+1)\times\ldots\ldots(n-1)\times n}{(1-i)\times(2-i)\times\ldots\times-1\times1\times2\times\ldots\times n-i}, \qquad i=1,2,3\ldots n\\&=\frac{(-1)^{i-1}n!}{i!(n-i)!}\end{aligned}$$

$$\prod_{\substack{j=1\\ i\neq j}}^{n}\frac{j}{j-i}=(-1)^{i-1}\,{}^{n}C_i, \qquad i=1,2,3\ldots n \tag{4.11}$$

$$\begin{aligned}\rho_m(a)&=\sum_{i=1}^{n}\frac{1}{(a+ih-a)^m}\prod_{\substack{j=1\\ i\neq j}}^{n}\frac{a-a-jh}{a+ih-a-jh}=\frac{1}{h^m}\sum_{i=1}^{n}\frac{1}{i^m}\prod_{\substack{j=1\\ i\neq j}}^{n}\frac{j}{j-i}, \quad m=1,2,3\ldots n\\&=\frac{1}{h^m}\sum_{i=1}^{n}\frac{1}{i^m}(-1)^{i-1}\,{}^{n}C_i\\&=\frac{1}{h^m}U^{(m)}, \quad m=1,2,3\ldots n\end{aligned} \tag{4.12}$$

Where $U^{(m)}=\sum_{i=1}^{n}\frac{1}{i^{m}}(-1)^{i-1}\,{}^{n}C_{i}$ (4.13)

**Note 4.2.** Suppose the data are spaced equally as $a+ih$ $(i=-m,\ldots-2,-1,0,1,2,\ldots n)$ then the value of $\rho_{m}(a)$ is calculated as follow as

$$\prod_{\substack{j=-m\\ i\neq 0,-i\neq j}}^{n}\frac{-j}{i-j}=\frac{m\ldots(n+1)n\ldots(i+1)(i-1)\times\ldots 2\times 1\times -1\times -2\times\ldots\times -n}{(-i+m)\ldots\times(-i+n+1)(-i+n)\ldots 1\times -1\times\ldots(-i+2)(-i+1)(-i-1)(-i-2)\ldots(-i-n)}$$

$$=(-1)^{i-1}\frac{(i+1)(i+2)\ldots m\times n!}{(m-i)!(i+1)(i+2)\ldots(i+n)}=A_{-i}\,,\qquad i=1,2,3\ldots m \tag{4.14}$$

$$\prod_{\substack{j=-m\\ i\neq 0,i\neq j}}^{n}\frac{-j}{i-j}=\frac{m(m-1)\ldots 2\times 1\times -1\times -2\times\ldots\times -(i-1)\times -(i+1)\times\ldots -n}{(i+m)(i+m-1)\times\ldots\times(i+2)(i+1)(i-1)(i-2)\times\ldots\times 1\times -1\times -2\times\ldots(i-n)}$$

$$=(-1)^{i-1}\frac{m!(i+1)(i+2)\times\ldots n}{(i+1)(i+2)\ldots(i+m)(n-i)!}=A_{i}\,,\qquad i=1,2,3\ldots n \tag{4.15}$$

$$\rho_{r}(a)=\sum_{\substack{i=-m\\ i\neq 0}}^{n}\frac{1}{(a+ih-a)^{r}}\prod_{\substack{j=-m\\ i\neq j,i,j\neq 0}}^{n}\frac{a-a-jh}{a+ih-a-jh}$$

$$=\frac{1}{h^{r}}\left(\sum_{i=-m}^{-1}\frac{1}{i^{r}}\prod_{\substack{j=-m\\ i\neq j,i,j\neq 0}}^{n}\frac{j}{j-i}+\sum_{i=1}^{n}\frac{1}{i^{r}}\prod_{\substack{j=-m\\ i\neq j,i,j\neq 0}}^{n}\frac{j}{j-i}\right)$$

$$=\frac{1}{h^{r}}\left((-1)^{r}\sum_{i=1}^{m}\frac{1}{i^{r}}\prod_{\substack{j=-m\\ -i\neq j,i,j\neq 0}}^{n}\frac{j}{j-i}+\sum_{i=1}^{n}\frac{1}{i^{r}}\prod_{\substack{j=-m\\ i\neq j,i,j\neq 0}}^{n}\frac{j}{j-i}\right)$$

$$\rho_{r}(a)=\frac{1}{h^{r}}\left((-1)^{r}\sum_{i=1}^{m}\frac{1}{i^{r}}A_{-i}+\sum_{i=1}^{n}\frac{1}{i^{r}}A_{i}\right),\quad r=1,2,3\ldots n \tag{4.16}$$

$$=\frac{1}{h^{r}}W^{(r)},\quad r=1,2,3\ldots n \tag{4.17}$$

Where $W^{(r)}=(-1)^{r}\sum_{i=1}^{m}\frac{1}{i^{r}}A_{-i}+\sum_{i=1}^{n}\frac{1}{i^{r}}A_{i}\,,\quad r=1,2,3\ldots n$ (4.18)

**Note 4.3.** Suppose the data are spaced as $a\pm ih$ $(i=0,1,\ldots n)$ then the value of $\rho_{m}(a)$ is calculated as follow as

If $m=n$

$$A_{-i}=(-1)^{i-1}\frac{(i+1)(i+2)\ldots n\times n!}{(n-i)!(i+1)(i+2)\ldots(i+n)} \tag{4.19}$$

$$=(-1)^{i-1}\frac{(i+1)(i+2)\ldots n\times(n-i)!(n-i+1)\ldots\times(n-1)\times n}{(n-i)!(i+1)(i+2)\ldots n(n+1)(n+2)\ldots(i+n)}$$

$$A_{-i}=(-1)^{i-1}\frac{n(n-1)\ldots(n-i+1)}{(n+1)(n+2)\ldots(i+n)}\,,\qquad i=1,2,3\ldots n \tag{4.20}$$

Similarly,

$$A_i = (-1)^{i-1} \frac{n(n-1)\dots(n-i+1)}{(n+1)(n+2)\dots(n+i)}, \qquad i = 1,2,3\dots n \tag{4.21}$$

$$tf \ r \ is \ even, \ \rho_r(a) = \frac{2V^{(r)}}{h^r} \quad \text{and} \quad if \ r \ is \ odd, \ \rho_r(a) = 0 \tag{4.22}$$

$$\text{Where } V^{(r)} = \sum_{i=1}^{n} \frac{A_i}{i^r} \qquad r = 1,2,3\dots k \tag{4.23}$$

**Corollary 4.1.2.** *Let* $a, a+h, a+2h, \dots, a+nh$ *are the distinct numbers in the interval* $[p,q]$, *spaced equally with* $h \neq 0$, $k \in W$ *and if* $f \in C^{n+k}[p,q]$, *then there exists a number* $\xi \in \{a, a+h, a+2h, \dots a+nh\}$ *with*

$$\begin{aligned} &\frac{f^{(k)}(a)}{k!}h^k + \frac{f^{(k-1)}(a)}{k-1!}U^{(1)}h^{k-1} + \frac{f^{(k-2)}(a)}{k-2!}U^{(2)}h^{k-2} + \dots + \frac{f(a)}{0!}U^{(k)} \\ &\qquad = \sum_{i=1}^{n} \frac{f(a+ih)}{i^k}(-1)^{i-1}\,{}^nC_i + (-1)^n h^{n+k} n! \frac{f^{(k+n)}(\xi)}{(k+n)!} \end{aligned} \tag{4.24}$$

Where $U^{(m)} = \frac{{}^nC_1}{1^m} - \frac{{}^nC_2}{2^m} + \dots + (-1)^{n-1}\frac{{}^nC_n}{n^m}$, $m = 1,2,\dots k$

**Proof.**

Putting $x = a$, $x_i = a + ih$, $i = 1,2,\dots n$ Applying Equation (4.2), we get

$$\begin{aligned} &\frac{f^{(k)}(a)}{k!}\rho_0(a) + \frac{f^{(k-1)}(a)}{k-1!}\rho_1(a) + \dots + \frac{f(a)}{0!}\rho_k(a) \\ &\qquad = \frac{1}{h^k}\sum_{i=1}^{n} \frac{f(a+ih)}{i^k} \prod_{\substack{j=1 \\ i\neq j}}^{n} \frac{j}{j-i} + (-1)^n h^n \frac{f^{(k+n)}(\xi)}{(k+n)!} \end{aligned} \tag{4.25}$$

Using equations (4.11), (4.12) and (4.13) (Refer Note 4.1) in (4.25), we obtain (4.24)

**Theorem 4.2.** *Let* $a+ih$ $(i = -m, \dots -2, -1, 0, 1, 2, \dots, n)$ *are the distinct numbers spaced equally in the interval* $[p,q]$, $h \neq 0$, $k \in W$ *and* $f \in C^{n+m+k}[p,q]$, *then for some* $\xi \in \{a-mh, \dots a-2h, , a-h, a, a+h, a+2h, \dots a+nh\}$

$$\begin{aligned} &\frac{f^{(k)}(a)}{k!}W^{(0)}h^k + \frac{f^{(k-1)}(a)}{k-1!}W^{(1)}h^{k-1} + \frac{f^{(k-2)}(a)}{k-2!}W^{(2)}h^{k-2} + \dots + \frac{f(a)}{0!}W^{(k)} \\ &\qquad = \sum_{i=1}^{m} \frac{(-1)^k A_{-i} f(a-ih)}{i^k} + \sum_{i=1}^{n} \frac{A_i f(a+ih)}{i^k} + (-1)^n h^{n+m+k} n! m! \frac{f^{(k+n+m)}(\xi)}{(k+n+m)!} \end{aligned} \tag{4.26}$$

Where $W^{(0)} = 1, W^{(r)} = (-1)^r \sum_{i=1}^{m} \frac{A_{-i}}{i^r} + \sum_{i=1}^{n} \frac{A_i}{i^r}$, $r = 1,2,\dots k$

$A_{-i} = \frac{(-1)^{i-1}(i+1)(i+2)\dots m.n!}{(m-i)!.(i+1)(i+2)\dots(i+n)}$, $i = 1,2,\dots m$ and $A_i = \frac{(-1)^{i-1} m!.(i+1)(i+2)\dots n}{(i+1)(i+2)\dots(i+m).(n-i)!}$, $i = 1,2,\dots n$

**Proof.**

Putting $a+ih$ $(i=-m\ldots,-2,-1,0,1,2,\ldots n)$, in applying Equation (4.2),

$$\frac{f^{(k)}(a)}{k!}\rho_0(a)+\frac{f^{(k-1)}(a)}{k-1!}\rho_1(a)+\ldots+\frac{f(a)}{0!}\rho_k(a)$$
$$=\frac{1}{h^k}\sum_{\substack{i=-m\\ i\neq 0}}^{n}\frac{f(a+ih)}{i^k}\prod_{\substack{j=-m\\ i\neq j,i,j\neq 0}}^{n}\frac{j}{j-i}+(-1)^n h^{n+m}\frac{f^{(k+n+m)}(\xi)n!m!}{(k+n+m)!} \tag{4.27}$$

Using equations (4.16), (4.17), and (4.18) (Refer Note 4.2) in (4.27) and simplifying, we get

$$\frac{f^{(k)}(a)}{k!}W^{(0)}h^k+\frac{f^{(k-1)}(a)}{k-1!}W^{(1)}h^{k-1}+\frac{f^{(k-2)}(a)}{k-2!}W^{(2)}h^{k-2}+\ldots+\frac{f(a)}{0!}W^{(k)}$$
$$=\sum_{i=1}^{n}\frac{f(a+ih)}{i^k}\prod_{\substack{j=-m\\ i\neq j,j\neq 0}}^{n}\frac{j}{j-i}+(-1)^k\sum_{i=1}^{m}\frac{f(a-ih)}{i^k}\prod_{\substack{j=-m\\ i\neq j,j\neq 0}}^{n}\frac{j}{j-i}+(-1)^n h^{n+m+k}\frac{f^{(n+m+k)}(\xi)n!m!}{(n+m+k)!} \tag{4.28}$$

Using equation (4.14) and (4.15) (Refer Note 4.2) in equation (4.28), we get

$$=\sum_{i=1}^{n}\frac{f(a+ih)A_i}{i^k}+(-1)^k\sum_{i=1}^{m}\frac{f(a-ih)A_{-i}}{i^k}+(-1)^n h^{n+m+k}\frac{f^{(n+m+k)}(\xi)n!m!}{(n+m+k)!} \tag{4.29}$$

Hence the proof.

**Corollary 4.2.1.** *Let* $a\pm ih$ $(i=0,1,2,\ldots n)$ *are the distinct numbers spaced equally in the interval* $[p,q]$, $h\neq 0$, $k\in W$ *and* $f\in C^{2n+k}[p,q]$, *then there exists a number* $\xi\in\{a-nh,\ldots a-2h,,a-h,a,a+h,a+2h,\ldots a+nh\}$ *with*

$$\frac{f^{(k)}(a)}{k!}h^k+2\left[\frac{f^{(k-2)}(a)}{(k-2)!}V^{(2)}h^{k-2}+\frac{f^{(k-4)}(a)}{(k-4)!}V^{(4)}h^{k-4}\ldots+\frac{f^{\psi}(a)}{\psi!}V^{(k)}\right]$$
$$=\sum_{i=1}^{n}\frac{f(a+ih)+(-1)^{\psi}f(a-ih)}{i^k}A_i+\frac{(-1)^n h^{2n+k}f^{(2n+k)}(\xi)n!n!}{(2n+k)!} \tag{4.30}$$

Where $\psi=\begin{cases}0, k\ is\ even\\ 1, k\ is\ odd\end{cases}$, $A_i=\frac{(-1)^{i-1}n(n-1)\ldots(n-i+1)}{(n+1)(n+2)\ldots(n+i)}$, $i=0,1,2,\ldots n$ and $V^{(m)}=\sum_{i=1}^{n}\frac{A_i}{i^m}$, $m=2,4\ldots$

**Proof.**

Substituting $m=n$ in equation (4.26), Using equations (4.20), (4.21), (4.22) and (4.23) (Refer Note 4.3)

After simplification, we obtain (4.30).

**Corollary 4.2.2.** *Let* $a\pm ih$ $(i=0,1,2,\ldots\infty)$ *are the distinct numbers spaced equally,* $h\neq 0$, $|h|<1$

*If* $f$ *is an infinitely differentiable function and* $\underset{n\to\infty}{Lt}\frac{f^{(2n+k)}(\xi)n!n!}{(2n+k)!}$ *exists for some* $\xi\in(-\infty,\infty)$ *then*

$$
\frac{f^{(k)}(a)}{k!}h^k + 2\left(\frac{f^{(k-2)}(a)}{(k-2)!}\sigma_2 h^{k-2} + \frac{f^{(k-4)}(a)}{(k-4)!}\sigma_4 h^{k-4} + \ldots + \frac{f^{\psi}(a)}{\psi!}\sigma_k\right)
$$

$$
= \sum_{i=1}^{\infty}(-1)^{i-1}\,\frac{f(a+ih)+(-1)^{\psi} f(a-ih)}{i^k} \tag{4.31}
$$

Where $\sigma_m = \frac{1}{1^m} - \frac{1}{2^m} + \frac{1}{3^m} - \ldots\infty$, $m = 2,4\ldots$ and $\psi = \begin{cases} 0, \, k \text{ is even} \\ 1, \, k \text{ is odd} \end{cases}$

**Proof .**

Letting $n \to \infty$ in Equation (4.30), Then

$$
\underset{n\to\infty}{Lt}\, A_i = \frac{(-1)^{i-1} n(n-1)\ldots(n-i+1)}{(n+1)(n+2)\ldots(n+i)} = (-1)^{i-1}\,\frac{\underset{n\to\infty}{Lt}\left(1-\frac{1}{n}\right)\left(1-\frac{2}{n}\right)\ldots\left(1-\frac{i-1}{n}\right)}{\underset{n\to\infty}{Lt}\left(1+\frac{1}{n}\right)\left(1+\frac{2}{n}\right)\left(1+\frac{3}{n}\right)\ldots\left(1+\frac{i}{n}\right)} = (-1)^{i-1}, \quad i = 1,2,\ldots n
$$

$\underset{n\to\infty}{Lt}\,\frac{f^{(2n+k)}(\xi)n!n!}{(2n+k)!}$ exists *and* $|h| < 1$, so that $\underset{n\to\infty}{Lt}\,(-1)^n h^{2n+k}\,\frac{f^{(2n+k)}(\xi)n!n!}{(2n+k)!} = 0$

Substituting all $A$'s in (4.30), we obtain equation (4.31).

**Theorem 4.3.** *Let* $x, x_0, x_1, \ldots x_n$ *are distinct numbers in the interval* $[p,q]$, $k \in W$ *and* $f \in C^{k+1}[p,q]$, *then*

$$
\frac{f^{(k)}(x)}{k!}\rho_0(x) + \frac{f^{(k-1)}(x)}{k-1!}\rho_1(x) + \ldots + \frac{f(x)}{0!}\rho_k(x)
$$

$$
+\frac{1}{k!}\sum_{i=0}^{n}\frac{1}{(x_i-x)^k}\prod_{\substack{j=0\\ i\neq j}}^{n}\frac{x-x_j}{x_i-x_j}\int_x^{x_i} f^{(k+1)}(t)(x_i-t)^k\,dt = \sum_{i=0}^{n}\frac{f(x_i)}{(x_i-x)^k}\prod_{\substack{j=0\\ i\neq j}}^{n}\frac{x-x_j}{x_i-x_j} \tag{4.32}
$$

Where $\rho_0(x) = 1$, $\rho_m(x) = \sum_{i=0}^{n}\frac{1}{(x_i-x)^m}\prod_{\substack{j=0\\ i\neq j}}^{n}\frac{x-x_j}{x_i-x_j}$, for $m = 1,2,3,\ldots k$

**Proof.**

Expanding $f(x_i)$ using Taylor formula

$$
f(x_i) = f(x) + \frac{f'(x)}{1!}(x_i - x) + \ldots + \frac{f^{(k)}(x)}{k!}(x_i-x)^k + \frac{1}{k!}\int_x^{x_i} f^{(k+1)}(t)(x_i-t)^k\,dt \tag{4.33}
$$

Dividing by $(x_i - x)^k$, we get

$$
\frac{f(x)}{(x_i-x)^k} + \frac{f'(x)}{(x_i-x)^{k-1}1!} + \ldots + \frac{f^{(k)}(x)}{k!} + \frac{1}{k!(x_i-x)^k}\int_x^{x_i} f^{(k+1)}(t)(x_i-t)^k\,dt = \frac{f(x_i)}{(x_i-x)^k} \tag{4.34}
$$

Multiplying by $\prod_{\substack{j=0\\ i\neq j}}^{n} \frac{x-x_j}{x_i-x_j}$ on both sides

$$\frac{f(x)}{(x_i-x)^k}\prod_{\substack{j=0\\ i\neq j}}^{n}\frac{x-x_j}{x_i-x_j}+\frac{f'(x)}{(x_i-x)^{k-1}1!}\prod_{\substack{j=0\\ i\neq j}}^{n}\frac{x-x_j}{x_i-x_j}+\ldots+\frac{f^{(k)}(x)}{k!}\prod_{\substack{j=0\\ i\neq j}}^{n}\frac{x-x_j}{x_i-x_j}$$
$$+\frac{1}{k!(x_i-x)^k}\prod_{\substack{j=0\\ i\neq j}}^{n}\frac{x-x_j}{x_i-x_j}\int_{x}^{x_i} f^{(k+1)}(t)(x_i-t)^k\,dt=\frac{f(x_i)}{(x_i-x)^k}\prod_{\substack{j=0\\ i\neq j}}^{n}\frac{x-x_j}{x_i-x_j} \quad (4.35)$$

Taking summation on both sides by taking $i=0,1,2,3,\ldots n$, then we obtain equation (4.32).

**Corollary 4.3.1.** *Let $x, x_0, x_1, \ldots x_n$ are distinct numbers in the interval $[p,q]$, $k\in W$ and $f\in C^{k+1}[p,q]$, then*

$$\frac{f^{(k)}(x)}{k!}\prod_{i=0}^{n}(x_i-x)^k\phi_0(x)+\prod_{i=0}^{n}(x_i-x)^{k-1}\frac{f^{(k-1)}(x)}{k-1!}\phi_1(x)+\ldots+\frac{f(x)}{0!}\phi_k(x)$$
$$+\frac{1}{k!}\sum_{i=0}^{n}\prod_{\substack{j=0\\ i\neq j}}^{n}\frac{(x_j-x)^{k+1}}{x_i-x_j}\int_{x}^{x_i} f^{(k+1)}(t)(x_i-t)^k\,dt=\sum_{i=0}^{n}f(x_i)\prod_{\substack{j=0\\ i\neq j}}^{n}\frac{(x_j-x)^{k+1}}{x_i-x_j} \quad (4.36)$$

Where $\phi_m(x)=\sum_{i=0}^{n}\prod_{\substack{j=0\\ i\neq j}}^{n}\frac{(x_j-x)^{m+1}}{x_i-x_j}$, $m=0,1,2,3,\ldots k$

**Corollary 4.3.2.** *Let $a, a+h, a+2h, \ldots, a+nh$ are the distinct numbers in the interval $[p,q]$, spaced equally with $h\neq 0$, $k\in W$ and if $f\in C^{k+1}[p,q]$, then*

$$\frac{f^{(k)}(a)}{k!}h^k+\frac{f^{(k-1)}(a)}{k-1!}h^{k-1}U^{(1)}+\frac{f^{(k-2)}(a)}{k-2!}h^2U^{(2)}+\ldots+\frac{f(a)}{0!}U^{(k)}$$
$$+\frac{1}{k!}\sum_{i=1}^{n}\frac{(-1)^{i-1}}{i^k}{}^nC_i\int_{a}^{a+ih} f^{(k+1)}(t)(a+ih-t)^k\,dt=\sum_{i=1}^{n}\frac{f(a+ih)}{i^k}(-1)^{i-1}\,{}^nC_i \quad (4.37)$$

*Where* $U^{(m)}=\frac{{}^nC_1}{1^k}-\frac{{}^nC_2}{2^k}+\ldots+(-1)^{n-1}\frac{{}^nC_n}{n^k}$, $\quad m=1,2,3,\ldots k$

**Corollary 4.3.3.** *Let $a+ih$ $(i=-m,\ldots-2,-1,0,1,2,\ldots n)$ are the distinct numbers spaced equally in the interval $[p,q]$, $h\neq 0$, $k\in W$ and $f\in C^{k+1}[p,q]$, then*

Let $W^{(r)}=\sum_{i=1}^{m}\frac{(-1)^r A_{-i}}{i^r}+\sum_{i=1}^{n}\frac{A_i}{i^r}$, $r=1,2,\ldots k$

$A_{-i}=\frac{(-1)^{i-1}(i+1)(i+2)\ldots m.n!}{(m-i)!(i+1)(i+2)\ldots(i+n)}$, $i=1,2,\ldots m$ *and* $A_i=\frac{(-1)^{i-1}m!(i+1)(i+2)\ldots n}{(i+1)(i+2)\ldots(i+m)(n-i)!}$, $i=1,2,\ldots n$ *then*

$$\frac{f^{(k)}(a)}{k!}h^k+\frac{f^{(k-1)}(a)}{k-1!}W^{(1)}h^{k-1}+\frac{f^{(k-2)}(a)}{k-2!}W^{(2)}h^{k-2}+\ldots+\frac{f(a)}{0!}W^{(k)}$$
$$+\frac{1}{k!}\left(\sum_{i=1}^{m}\frac{(-1)^k A_{-i}}{i^k}\int_a^{a-ih}f^{(k+1)}(t)(a-ih-t)^k\,dt+\sum_{i=1}^{n}\frac{A_i}{i^k}\int_a^{a+ih}f^{(k+1)}(t)(a+ih-t)^k\,dt\right) \tag{4.38}$$
$$=\sum_{i=1}^{m}\frac{(-1)^k A_{-i}f(a-ih)}{i^k}+\sum_{i=1}^{n}\frac{A_i f(a+ih)}{i^k}$$

**Corollary 4.3.4.** *Let* $a\pm ih$ $(i=0,1,2,\ldots n)$ *are the distinct numbers spaced equally in the interval* $[p,q]$, $h\neq 0$, $k\in W$ *and* $f\in C^{k+1}[p,q]$,

$$\frac{f^{(k)}(a)}{k!}h^k+2\left[\frac{f^{(k-2)}(a)}{(k-2)!}V^{(2)}h^{k-2}+\frac{f^{(k-4)}(a)}{(k-4)!}V^{(4)}h^{k-4}+\ldots+\frac{f^{\psi}(a)}{\psi!}V^{(k)}\right]$$
$$+\frac{1}{k!}\sum_{i=1}^{n}\frac{A_i}{i^k}\left(\int_a^{a+ih}f^{(k+1)}(t)(a+ih-t)^k\,dt+(-1)^{\psi}\int_a^{a-ih}f^{(k+1)}(t)(a-ih-t)^k\,dt\right) \tag{4.39}$$
$$=\sum_{i=1}^{n}A_i\,\frac{f(a+ih)+(-1)^{\psi}f(a-ih)}{i^k}$$

Where $\psi=\begin{cases}0, & k \text{ is even}\\ 1, & k \text{ is odd}\end{cases}$ , $V^{(m)}=\sum_{i=1}^{n}\frac{A_i}{i^m}$, $m=2,4\ldots$ and $A_i=\frac{(-1)^{i-1}n(n-1)\ldots(n-i+1)}{(n+1)(n+2)\ldots(n+i)}$, $i=1,2,3,\ldots n$

**Corollary 4.3.5.** *Let* $a\pm ih$ $(i=0,1,2,\ldots\infty)$ *are the distinct numbers spaced equally,* $h\neq 0$, $|h|<1$
*If* $f$ *is an infinitely differentiable function , then*

$$\frac{f^{(k)}(a)}{k!}h^k+2\left[\frac{f^{(k-2)}(a)}{(k-2)!}\sigma_2 h^{k-2}+\frac{f^{(k-4)}(a)}{(k-4)!}\sigma_4 h^{k-4}+\ldots+\frac{f^{\psi}(a)}{\psi!}\sigma_k\right]$$
$$+\frac{1}{k!}\sum_{i=1}^{\infty}\frac{(-1)^{i-1}}{i^k}\left(\int_a^{a+ih}f^{(k+1)}(t)(a+ih-t)^k\,dt+(-1)^{\psi}\int_a^{a-ih}f^{(k+1)}(t)(a-ih-t)^k\,dt\right) \tag{4.40}$$
$$=\sum_{i=1}^{\infty}(-1)^{i-1}\,\frac{f(a+ih)+(-1)^{\psi}f(a-ih)}{i^k}$$

Where, $\sigma_m=\frac{1}{1^m}-\frac{1}{2^m}+\frac{1}{3^m}-\ldots\infty$ , $m=2,4\ldots$ and $\psi=\begin{cases}0, & k \text{ is even}\\ 1, & k \text{ is odd}\end{cases}$

**Corollary 4.3.6.** *Let* $x, x_0, x_1,\ldots x_n$ *are distinct numbers in the interval* $[p,q]$, $k\in W$ *and* $f\in C^{n+k+1}[p,q]$, *then for some* $\xi\in\{x, x_0, x_1,\ldots x_n\}$

(i). $$f[\underbrace{x,\ldots x}_{k+1\ times},x_0,x_1,x_2\ldots x_n]\prod_{i=0}^{n}(x-x_i)=-\frac{1}{k!}\sum_{i=0}^{n}\frac{1}{(x_i-x)^k}\prod_{\substack{j=0\\ i\neq j}}^{n}\frac{x-x_j}{x_i-x_j}\int_x^{x_i}f^{(k+1)}(t)(x_i-t)^k\,dt \tag{4.41}$$

(ii). $$\frac{f^{(n+k+1)}(\xi)}{(n+k+1)!}\prod_{i=0}^{n}(x-x_i) = -\frac{1}{k!}\sum_{i=0}^{n}\frac{1}{(x_i-x)^k}\prod_{\substack{j=0\\ i\neq j}}^{n}\frac{x-x_j}{x_i-x_j}\int_{x}^{x_i} f^{(k+1)}(t)(x_i-t)^k\,dt \tag{4.42}$$

**Corollary 4.3.7.** *Let* $a \pm ih$ $(i = 0,1,2,\ldots\infty)$ *are the distinct numbers spaced equally,* $h \neq 0$, $|h| < 1$, $k \in W$

*If* $f$ *is an infinitely differentiable function and* $\underset{n\to\infty}{Lt}\ \frac{f^{(2n+k)}(\xi)n!n!}{(2n+k)!}$ *extis for some* $\xi \in (-\infty,\infty)$, *then*

$$\sum_{i=1}^{\infty}\frac{(-1)^{i-1}}{i^k}\int_{a}^{a+ih} f^{(k+1)}(t)(a+ih-t)^k\,dt = (-1)^{\psi}\sum_{i=1}^{\infty}\frac{(-1)^{i-1}}{i^k}\int_{a-ih}^{a} f^{(k+1)}(t)(a-ih-t)^k\,dt \tag{4.43}$$

**Proof.**

Equating (4.31) and (4.40) and after simplification, we obtain (4.43).

If $n = 0$, and $k = 0$ we get Taylor's formula and bary centric interpolation formula respectively. We get higher order differentiation formula recursively, for different values of *k*, thus the recursive approach is more general formula for interpolation and numerical differentiation.

**4.2. Numerical differentiation formulas using linear combination of divided differences.**

**Theorem 4.4.** *Let* $x_0, x_1, \ldots x_n$ *and* $x$ *are distinct numbers in the interval* $[p,q]$, $k \in W$ *and* $f \in C^{n+k+1}[p,q]$, *then*

(i). *If functional value at x is known*

$$\begin{aligned}\frac{f^{(k)}(x)}{k!} = &\sum_{i_1<i_2<\ldots<i_k=0}^{n} f[x_{i_1}, x_{i_2}\ldots x_{i_k}, x]\prod_{\substack{j=0\\ i_1<i_2<\ldots<i_k\neq j}}^{n}\frac{(x-x_j)^k}{(x_{i_1}-x_j)(x_{i_2}-x_j)\cdots(x_{i_k}-x_j)}\\ &+ f[x,x,x_0,x_1,\ldots x_n]\frac{P^{(k-1)}(x)}{k-1!} + \ldots + f[\underbrace{x,\ldots x}_{k\ times}, x_0,x_1,\ldots x_n]\frac{P^{(1)}(x)}{1!}\\ &+ f[\underbrace{x,\ldots,x}_{k+1\ times}, x_0,x_1,\ldots x_n]P(x)\end{aligned} \tag{4.44}$$

(ii). *If functional value at x is not known*

$$\begin{aligned}\frac{f^{(k)}(x)}{k!} = &\sum_{i_1<i_2<\ldots<i_{k+1}=0}^{n} f[x_{i_1}, x_{i_2}\ldots x_{i_k}, x_{i_{k+1}}]\prod_{\substack{j=0\\ i_1<i_2<\ldots<i_{k+1}\neq j}}^{n}\frac{(x-x_j)^{k+1}}{(x_{i_1}-x_j)(x_{i_2}-x_j)\cdots(x_{i_{k+1}}-x_j)}\\ &+ f[x,x_0,x_1,\ldots x_n]\frac{P^{(k)}(x)}{k!} + \ldots + f[\underbrace{x,\ldots,x}_{k\ times}, x_0,x_1,\ldots x_n]\frac{P^{(1)}(x)}{1!}\\ &+ f[\underbrace{x,\ldots x}_{k+1\ times}, x_0,x_1,\ldots x_n]P(x)\end{aligned} \tag{4.45}$$

Where $P(x) = \prod_{i=0}^{n}(x-x_i)$ and $P^{(m)}(x)$ is $m^{th}$ order differentiation of $P(x)$, $m = 1,2,3,\ldots k$

**Proof.**

Expanding $f[\underbrace{x,\ldots,x}_{k+1\ times}]$ using Theorem 2.1

$$f[\underbrace{x,\ldots,x}_{k+1\ times}]=\sum_{i_1=0}^{n} f[x_{i_1},\underbrace{x,\ldots x}_{k\ times}]\prod_{\substack{j_1=0\\ i_1\neq j_1}}^{n}\frac{x-x_{j_1}}{x_{i_1}-x_{j_1}}+f[\underbrace{x,\ldots,x}_{k+1\ times},x_0,x_1,\ldots x_n]\prod_{j_1=0}^{n}(x-x_{j_1}) \tag{4.46}$$

Again applying Theorem 2.1 for $f[x_{i_1},\underbrace{x,\ldots,x}_{k\ times}]$ and after simplification

$$\frac{f^{(k)}(x)}{k!}=\sum_{i_1=0}^{n}\sum_{\substack{i_2=0\\ i_1\neq i_2}}^{n} f[x_{i_1},x_{i_2},\underbrace{x\ldots,x}_{k-1\ times}]\prod_{\substack{j_2=0\\ i_2,i_1\neq j_2}}^{n}\frac{x-x_{j_2}}{x_{i_2}-x_{j_2}}\prod_{\substack{j_1=0\\ i_1\neq j_1}}^{n}\frac{x-x_{j_1}}{x_{i_1}-x_{j_1}}$$

$$+f[\underbrace{x,\ldots x}_{k\ times},x_0,x_1,\ldots x_n]\sum_{i_1=0}^{n}\prod_{\substack{j_2=0\\ i_1\neq j_2}}^{n}x-x_{j_2}\prod_{\substack{j_1=0\\ i_1\neq j_1}}^{n}\frac{x-x_{j_1}}{x_{i_1}-x_{j_1}}+f[\underbrace{x,\ldots x}_{k+1\ times},x_0,x_1,\ldots x_n]\prod_{j_1=0}^{n}(x-x_{j_1})$$

Repeating this, k times

$$=\sum_{i_1=0}^{n}\sum_{\substack{i_2=0\\ i_1\neq i_2}}^{n}\cdots\sum_{\substack{i_k=0\\ i_1\neq i_2\neq\ldots\neq i_k}}^{n}[x_{i_1},x_{i_2},\ldots,x_{i_k},x]\prod_{\substack{j_k=0\\ i_1,i_2,\ldots,i_k\neq j_k}}^{n}\frac{x-x_{j_k}}{x_{i_k}-x_{j_k}}\cdots\prod_{\substack{j_2=0\\ i_2,i_1\neq j_2}}^{n}\frac{x-x_{j_2}}{x_{i_2}-x_{j_2}}\prod_{\substack{j_1=0\\ i_1\neq j_1}}^{n}\frac{x-x_{j_1}}{x_{i_1}-x_{j_1}}$$

$$+f[x,x,x_0,x_1,\ldots x_n]\sum_{i_1=0}^{n}\sum_{\substack{i_2=0\\ i_1\neq i_2}}^{n}\cdots\sum_{\substack{i_{k-1}=0\\ i_1\neq i_2,\neq\ldots\neq i_{k-1}}}^{n}\prod_{\substack{j_k=0\\ i_1,i_2,\ldots,i_{k-1}\neq j_k}}^{n}x-x_{j_k}\prod_{\substack{j_{k-1}=0\\ i_1,i_2,\ldots,i_{k-1}\neq j_{k-1}}}^{n}\frac{x-x_{j_{k-1}}}{x_{i_{k-1}}-x_{j_{k-1}}}\cdots\prod_{\substack{j_1=0\\ i_1\neq j_1}}^{n}\frac{x-x_{j_1}}{x_{i_1}-x_{j_1}} \tag{4.47}$$

$$+\ldots+f[\underbrace{x,\ldots x}_{k\ times},x_0,x_1,\ldots x_n]\sum_{i_1=0}^{n}\prod_{\substack{j_2=0\\ i_1\neq j_2}}^{n}x-x_{j_2}\prod_{\substack{j_1=0\\ i_1\neq j_1}}^{n}\frac{x-x_{j_1}}{x_{i_1}-x_{j_1}}+f[\underbrace{x,\ldots x}_{k+1\ times},x_0,x_1,\ldots x_n]\prod_{j_1=0}^{n}(x-x_{j_1})$$

For any permutation $(i_1,i_2,\ldots,i_n)$ of $(j_1,j_2,\ldots,j_n)$, all $j$'s $=0,1,2,3,\ldots n$, We have

$$f[x_{i_1},x_{i_2}\ldots x_{i_k},x]=f[x_{j_1},x_{j_2}\ldots x_{j_k},x] \text{ and } \sum_{p=1}^{m}\prod_{q=1,p\neq q}^{m}\frac{x-x_{i_q}}{x_{i_p}-x_{i_q}}=1,\ m\in N$$

Using above equations and after simplification, we get

$$\frac{f^{(k)}(x)}{k!}=\sum_{i_1=0}^{n}\sum_{\substack{i_2=0\\ i_1<i_2}}^{n}\cdots\sum_{\substack{i_k=0\\ i_1<i_2<\ldots<i_k}}^{n} f[x_{i_1},x_{i_2}\ldots x_{i_k},x]\prod_{\substack{j=0\\ i_1<i_2<\ldots i_k\neq j}}^{n}\frac{(x-x_j)^k}{(x_{i_1}-x_j)(x_{i_2}-x_j)\cdots(x_{i_k}-x_j)}$$

$$+f\left[x,x,x_0,x_1,\ldots x_n\right]\sum_{i_1=0}^{n}\sum_{\substack{i_2=0\\ i_1<i_2}}^{n}\cdots\sum_{\substack{i_{k-1}=0\\ i_1<i_2<\ldots i_{k-1}}}^{n}\prod_{\substack{j=0\\ i_1<i_2<\ldots i_{k-1}\neq j}}^{n}\frac{(x-x_j)^k}{(x_{i_1}-x_j)(x_{i_2}-x_j)\ldots(x_{i_{k-1}}-x_j)} \tag{4.48}$$

$$+\ldots+f[\underbrace{x,\ldots,x}_{k\ times},x_0,x_1,\ldots x_n]\sum_{i_1=0}^{n}\prod_{\substack{j=0\\ i_1\neq j}}^{n}\frac{(x-x_j)^2}{(x_{i_1}-x_j)}+f[\underbrace{x,\ldots,x}_{k+1\ times},x_0,x_1,\ldots x_n]\prod_{j_1=0}^{n}(x-x_{j_1})$$

again expanding (4.47) using Theorem 2.1 after simplification we get

$$\frac{f^{(k)}(x)}{k!}=\sum_{i_1=0}^{n}\sum_{\substack{i_2=0\\i_1<i_2}}^{n}\cdots\sum_{\substack{i_{k+1}=0\\i_1<i_2<\ldots<i_{k+1}}}^{n} f[x_{i_1},x_{i_2}\ldots x_{i_k},x_{i_{k+1}}]\prod_{\substack{j=0\\i_1<i_2<\ldots i_{k+1}\neq j}}^{n}\frac{(x-x_j)^{k+1}}{(x_{i_1}-x_j)(x_{i_2}-x_j)\cdots(x_{i_{k+1}}-x_j)}$$

$$+f[x,x_0,x_1,\ldots x_n]\sum_{i_1=0}^{n}\sum_{\substack{i_2=0\\i_1<i_2}}^{n}\cdots\sum_{\substack{i_k=0\\i_1<i_2<\ldots i_k}}^{n}\prod_{\substack{j=0\\i_1<i_2\ldots i_k\neq j}}^{n}\frac{(x-x_j)^{k+1}}{(x_{i_1}-x_j)(x_{i_2}-x_j)\ldots(x_{i_k}-x_j)} \tag{4.49}$$

$$+\ldots+f[\underbrace{x,\ldots,x}_{k\ times},x_0,x_1,\ldots x_n]\sum_{i_1=0}^{n}\prod_{\substack{j=0\\i_1\neq j}}^{n}\frac{(x-x_j)^2}{(x_{i_1}-x_j)}+f[\underbrace{x,\ldots,x}_{k+1\ times},x_0,x_1,\ldots x_n]\prod_{j_1=0}^{n}(x-x_{j_1})$$

Now, replace $f(x)$ by $P(x)=\prod_{i=0}^{n}(x-x_i)$, then $P(x_i)=0$, Substituting in (4.49) we get

$$\frac{P^{(m)}(x)}{m!}=\sum_{i_1=0}^{n}\sum_{\substack{i_2=0\\i_1<i_2}}^{n}\cdots\sum_{\substack{i_m=0\\i_1<i_2<\ldots<i_m}}^{n}\prod_{\substack{j=0\\i_1,i_2,\ldots,i_m\neq j}}^{n}\frac{(x-x_j)^{m+1}}{(x_{i_1}-x_j)(x_{i_2}-x_j)\ldots(x_{i_m}-x_j)},\quad m=0,1,2,3,\ldots,k \tag{4.50}$$

Applying (4.50) for $m=0,1,2,3,\ldots k-1$ in (4.48) we get,

$$=\sum_{i_1=0}^{n}\sum_{\substack{i_2=0\\i_1<i_2}}^{n}\cdots\sum_{\substack{i_k=0\\i_1<i_2<\ldots<i_k}}^{n} f[x_{i_1},x_{i_2}\ldots x_{i_k},x]\prod_{\substack{j=0\\i_1,i_2,\ldots i_k\neq j}}^{n}\frac{(x-x_j)^{k}}{(x_{i_1}-x_j)(x_{i_2}-x_j)\cdots(x_{i_k}-x_j)}$$

$$+f[x,x,x_0,x_1,\ldots x_n]\frac{P^{(k-1)}}{k-1!}(x)+\ldots+f[\underbrace{x,\ldots x}_{k\ times},x_0,x_1,\ldots x_n]\frac{P^{(1)}(x)}{1!}+f[\underbrace{x,\ldots x}_{k+1\ times},x_0,x_1,\ldots x_n]P(x) \tag{4.51}$$

Using the notation $\sum_{i_1=0}^{n}\sum_{\substack{i_2=0\\i_1<i_2}}^{n}\cdots\sum_{\substack{i_k=0\\i_1<i_2<\ldots<i_k}}^{n} f[x_{i_1},x_{i_2}\ldots x_{i_k},x]=\sum_{i_1<i_2<\ldots<i_k=0}^{n} f[x_{i_1},x_{i_2}\ldots x_{i_k},x]$ in (4.51), we get (4.44)

Similarly, from (4.49), we can easily obtain (4.45).

**Corollary 4.4.1.** *Let* $x_0,x_1,\ldots x_n$ *and* $x$ *are distinct number in the interval* $[p,q]$, $k\in W$ *and* $f\in C^{n+k+1}[p,q]$,
(i). *If functional value at* $x$ *is known, for some* $\xi_i$, $\xi_i\in\{x,x_0,x_1,\ldots x_n\}$, $i=1,2,3,\ldots,k$

$$\frac{f^{(k)}(x)}{k!}=\sum_{i_1<i_2<\ldots<i_k=0}^{n} f[x_{i_1},x_{i_2}\ldots x_{i_k},x]\prod_{\substack{j=0\\i_1<i_2<\ldots<i_k\neq j}}^{n}\frac{(x-x_j)^{k}}{(x_{i_1}-x_j)(x_{i_2}-x_j)\cdots(x_{i_k}-x_j)}$$

$$+\frac{f^{(n+2)}(\xi_k)}{(n+2)!}\frac{P^{(k-1)}(x)}{k-1!}+\ldots+\frac{f^{(k+n)}(\xi_2)}{(k+n)!}\frac{P^{(1)}(x)}{1!}+\frac{f^{(k+n+1)}(\xi_1)}{(k+n+1)!}P(x) \tag{4.52}$$

(ii). *If functional value at* $x$ *is not known, for some* $\xi_i$, $\xi_i\in\{x,x_0,x_1,\ldots x_n\}$, $i=1,2,3,\ldots,k+1$

$$\frac{f^{(k)}(x)}{k!}=\sum_{i_1<i_2<\ldots<i_{k+1}=0}^{n} f[x_{i_1},x_{i_2}\ldots x_{i_k},x_{i_{k+1}}]\prod_{\substack{j=0\\i_1<i_2<\ldots<i_{k+1}\neq j}}^{n}\frac{(x-x_j)^{k+1}}{(x_{i_1}-x_j)(x_{i_2}-x_j)\cdots(x_{i_{k+1}}-x_j)}$$

$$+\frac{f^{(n+1)}(\xi_{k+1})}{(n+1)!}\frac{P^{(k)}(x)}{k!}+\ldots+\frac{f^{(k+n)}(\xi_2)}{(k+n)!}\frac{P^{(1)}(x)}{1!}+\frac{f^{(k+n+1)}(\xi_1)}{(k+n+1)!}P(x) \tag{4.53}$$

**Corollary 4.4.2.** *Let* $a, a+h, a+2h, \ldots, a+nh$ *are the distinct numbers in the interval* $[p,q]$, *spaced equally with* $h \neq 0$ $x = a+th \in [p,q]$, *and if* $f \in C^{n+k+1}[p,q]$, *with*

(i). *If functional value at* $x$ *is known, then*

$$\frac{f^{(k)}(x)}{k!} = \frac{1}{h^k} \sum_{i_1<i_2<\ldots<i_k=0}^{n} f_I[i_1, i_2, \ldots i_k, t] \prod_{\substack{j=0 \\ i_1<i_2<\ldots<i_k \neq j}}^{n} \frac{(t-j)^k}{(i_1-j)(i_2-j)\cdots(i_k-j)} + O(h^{n+1-k}) \tag{4.54}$$

(ii). *If functional value at x is not known, then*

$$\frac{f^{(k)}(x)}{k!} = \frac{1}{h^k} \sum_{i_1<i_2<\ldots<i_{k+1}=0}^{n} f_I[i_1, i_2, \ldots, i_{k+1}] \prod_{\substack{j=0 \\ i_1<i_2<\ldots<i_{k+1} \neq j}}^{n} \frac{(t-j)^{k+1}}{(i_1-j)(i_2-j)\cdots(i_{k+1}-j)} + O(h^{n+1-k}) \tag{4.55}$$

**Corollary 4.4.3.** *Let* $a, a+h, a+2h, \ldots, a+nh$ *are the distinct numbers in the interval* $[p,q]$, *spaced equally with* $h \neq 0$ *and if* $f \in C^{n+k}[p,q]$, *with*

$$\frac{f^{(k)}(a)}{k!} = \frac{1}{h^k} \sum_{i_1<i_2<\ldots<i_k=1}^{n} f_I[i_1, i_2, \ldots i_k, 0](-1)^{i_1+i_2+\ldots i_k-k} \frac{{}^nC_{i_1} {}^nC_{i_2} \ldots {}^nC_{i_k}}{(i_1 i_2 \ldots i_k)^{k-1}} \pi_1 \pi_2 \ldots \pi_k + O(h^{n-k}) \quad , \pi_z = \prod_{\substack{j=1 \\ z \neq j}}^{n} (i_z - i_j) \tag{4.56}$$

Suffix $I$ denotes the divided differences by integer by arguments

**Remark 4.1**

If $f(x) = x^k$, then from (4.44) and (4.56), we obtain following two results

(i).
$$\sum_{i_1<i_2<\ldots<i_k=0}^{n} \prod_{\substack{j=0 \\ i_1<i_2<\ldots<i_k \neq j}}^{n} \frac{(x-x_j)^k}{(x_{i_1}-x_j)(x_{i_2}-x_j)\cdots(x_{i_k}-x_j)} = 1 \tag{4.57}$$

(ii).
$$\sum_{i_1<i_2<\ldots<i_k=1}^{n} (-1)^{i_1+i_2+\ldots i_k-k} \frac{{}^nC_{i_1} {}^nC_{i_2} \ldots {}^nC_{i_k}}{(i_1 i_2 \ldots i_k)^{k-1}} \pi_1 \pi_2 \ldots \pi_k = 1, \text{ where } \pi_z = \prod_{\substack{j=1 \\ z \neq j}}^{n} (i_z - i_j) \tag{4.58}$$

## 5. Comparisons with other formulas of numerical differentiation

To compare the differentiation formulas based on recursive approach (discussed in the previous section) with former formulas of numerical differentiation. In this section, we derive them in terms of functional values at given nodes. While using recursive formula given in equation (4.1), to find differentiation of $t^{th}$ order. We get the following form for unevenly spaced points

$$\frac{f^{(t)}(x)}{t!} = M_t(x)a_0 + M_{t-1}(x)a_1 + \ldots + M_1(x)a_{t-1} + a_t f(x), \quad t = 1,2,\ldots \tag{5.1}$$

Where unknown $a's \in R$, $L_i(x) = \prod_{\substack{j=0 \\ i \neq j}}^{n} \frac{x-x_j}{x_i - x_j}$ for $i = 0,1,2,\ldots,n$

$$M_k(x)=\sum_{i=0}^{n}\frac{f(x_i)}{(x_i-x)^k}L_i(x)+f[\underbrace{x,\ldots x}_{k+1\ times},x_0,x_1,x_2\ldots x_n]\prod_{i=0}^{n}(x-x_i)\,,\ k=1,2,\ldots,t$$

For evenly spaced points,

$$\frac{f^{(t)}(x)}{t!}=\frac{1}{h^t}\left(\hat{M}_t(a)+\hat{M}_{t-1}(a)\hat{a}_1+\ldots+\hat{M}_1(a)\hat{a}_{t-1}+f(a)\hat{a}_t\right),\ t=1,2,\ldots \tag{5.2}$$

Where unknown $\hat{a}'s\in R$ and $\hat{M}_k(a)=\sum_{i=1}^{n}\frac{f(a+ih)}{i^k}(-1)^{i-1}\,{}^nC_i+(-1)^n h^{n+k}\frac{f^{(k+n)}(\xi_{t-k})}{(k+n)!}$, $k=1,2,\ldots,t$

**Theorem 5.1.** *Let* $x,x_0,x_1,\ldots x_n$ *are distinct numbers in the interval* $[p,q]$, $t\in W$ *and* $f\in C^{n+t+1}[p,q]$, *then*

(i)
$$\frac{f^{(t)}(x)}{t!}=\sum_{i=0}^{n}\left(\frac{a_0}{(x_i-x)^t}+\frac{a_1}{(x_i-x)^{t-1}}+\ldots\frac{a_{t-1}}{x_i-x}\right)f(x_i)L_i(x)+f(x)$$
$$+\prod_{i=0}^{n}(x-x_i)\left(a_0f[\underbrace{x,\ldots x}_{t+1\ times},x_0,x_1,\ldots x_n]+a_1f[\underbrace{x,\ldots x}_{t\ times},x_0,x_1,\ldots x_n]+\ldots+a_{t-1}f[x,x,x_0,x_1,\ldots x_n]\right) \tag{5.3}$$

(ii)
$$\frac{f^{(t)}(x)}{t!}=\sum_{i=0}^{n}\left(\frac{a_0}{(x_i-x)^t}+\frac{a_1}{(x_i-x)^{t-1}}+\ldots\frac{a_{t-1}}{x_i-x}\right)f(x_i)L_i(x)+f(x)$$
$$+\prod_{i=0}^{n}(x-x_i)\left(a_0\frac{f^{(t+n+1)}(\xi_0)}{(t+n+1)!}+a_1\frac{f^{(t+n)}(\xi_1)}{(t+n)!}+\ldots+a_{t-1}\frac{f^{(n+2)}(\xi_{t-1})}{(n+2)!}\right) \tag{5.4}$$

(iii)
$$\frac{f^{(t)}(x)}{t!}=\sum_{i=0}^{n}\left(\frac{a_0}{(x_i-x)^t}+\frac{a_1}{(x_i-x)^{t-1}}+\ldots\frac{a_{t-1}}{x_i-x}+a_t\right)f(x_i)L_i(x)$$
$$+\prod_{i=0}^{n}(x-x_i)\left(a_0\frac{f^{(t+n+1)}(\xi_0)}{(t+n+1)!}+a_1\frac{f^{(t+n)}(\xi_1)}{(t+n)!}+\ldots+a_t\frac{f^{(n+1)}(\xi_t)}{(n+1)!}\right) \tag{5.5}$$

Where $L_i(x)=\prod_{\substack{j=0\\ i\neq j}}^{n}\frac{x-x_j}{x_i-x_j}$ for $i=0,1,2,\ldots,n$, $\rho_k=\sum_{i=0}^{n}\frac{L_i(x)}{(x_i-x)^k}$ for $k=1,2,3,4,\ldots,t$

$a_0=1$, $a_k=-(\rho_ka_0+\rho_{k-1}a_1+\rho_{k-2}a_2+\ldots+\rho_1a_{k-1})$, $k=1,2,\ldots,t$ and $\xi_i\in\{x,x_0,x_1,\ldots x_n\}$, $i=0,1,2,\ldots,t$

**Proof:**

We can write for some $t=0,1,2,\ldots$ using equation (5.1)

$$\frac{f^{(t)}(x)}{t!}=M_t(x)a_0+M_{t-1}(x)a_1+\ldots+M_1(x)a_{t-1}+a_tf(x) \tag{5.6}$$

Where $M_k(x)=\sum_{i=0}^{n}\frac{f(x_i)}{(x_i-x)^k}L_i(x)+f[\underbrace{x,\ldots x}_{k+1\ times},x_0,x_1,x_2\ldots x_n]\prod_{i=0}^{n}(x-x_i)$, $k=1,2,\ldots,t$

$$\frac{f^{(t)}(x)}{t!} = \sum_{i=0}^{n}\left(\frac{a_0}{(x_i-x)^t}+\frac{a_1}{(x_i-x)^{t-1}}+\ldots+\frac{a_{t-1}}{x_i-x}\right)f(x_i)L_i(x)+a_t f(x)$$
$$+\prod_{i=0}^{n}(x-x_i)\left(a_0 f[\underbrace{x,\ldots x}_{t+1\ times},x_0,x_1,\ldots x_n]+a_1 f[\underbrace{x,\ldots x}_{t\ times},x_0,x_1,\ldots x_n]+\ldots+a_{t-1}f[x,x,x_0,x_1,\ldots x_n]\right) \quad (5.7)$$

To find $a$'s value, Put $f(x)=1$ in (5.7), If $t=0$, then $a_0=1$ and for $t=1,2,3,4,\ldots$ we have

$$\sum_{i=0}^{n}\left(\frac{a_0}{(x_i-x)^t}+\frac{a_1}{(x_i-x)^{t-1}}+\ldots\frac{a_{t-1}}{x_i-x}\right)L_i(x)+a_t=0$$

$$a_0\sum_{i=0}^{n}\frac{L_i(x)}{(x_i-x)^t}+a_1\sum_{i=0}^{n}\frac{L_i(x)}{(x_i-x)^{t-1}}+\ldots+a_{t-1}\sum_{i=0}^{n}\frac{L_i(x)}{x_i-x}+a_t=0,\ \text{Say}\ \sum_{i=0}^{n}\frac{L_i(x)}{(x_i-x)^m}=\rho_m,\ m=1,2,\ldots,t$$

$$a_0\rho_t+a_1\rho_{t-1}+\ldots+a_{t-1}\rho_1+a_t=0 \quad (5.8)$$

Using (5.8) and (5.7) we obtain (5.3).

We know that, for some $\xi_i\in\{x,x_0,x_1,\ldots x_n\}$, $i=0,1,2,\ldots,t$

$$f[\underbrace{x,\ldots x}_{i+1\ times},x_0,x_1,x_2\ldots x_n]=\frac{f^{(n+i+1)}(\xi_{t-i})}{(n+i+1)!} \quad (5.9)$$

Substituting (5.9) in (5.3) we obtain (5.4).

Again using Lagrange interpolation formula for $f(x)$, and simplifying we get (5.5)

**Corollary 5.1.1.** *Let* $a,a+h,a+2h,\ldots,a+nh$ *are the distinct numbers in the interval* $[p,q]$, *spaced equally with* $h\neq 0$, $t\in W$ *and if* $f\in C^{n+t}[p,q]$, *then for* $\xi_i\in\{a,a+h,a+2h,\ldots a+nh\}$, $i=0,1,2,\ldots,t-1$

$$\frac{f^{(t)}(a)}{t!}=\frac{1}{h^t}\left(\sum_{i=1}^{n}(-1)^{i-1}\,{}^{n}C_i f(a+ih)\left(\frac{\hat{a}_0}{i^t}+\frac{\hat{a}_1}{i^{t-1}}+\ldots\frac{\hat{a}_{t-1}}{i}\right)+\hat{a}_t f(a)\right.$$
$$\left.+(-1)^n h^n n!\left(\hat{a}_0\frac{f^{(t+n)}(\xi_0)}{(t+n)!}+\frac{\hat{a}_1}{h}\frac{f^{(t+n-1)}(\xi_1)}{(t+n-1)!}+\ldots+\frac{\hat{a}_{t-1}}{h^{t-1}}\frac{f^{(n+1)}(\xi_{t-1})}{n+1!}\right)\right) \quad t=0,1,2,\ldots \quad (5.10)$$

*Where* $U^{(k)}=\sum_{i=1}^{n}\frac{(-1)^i\,{}^nC_i}{i^k}$ , $k=1,2,\ldots,t$

$\hat{a}_0=1$, $\hat{a}_k=-(U^{(k)}\hat{a}_0+U^{(k-1)}\hat{a}_1+U^{(k-2)}\hat{a}_2+\ldots+U^{(1)}\hat{a}_{k-1})$, $k=1,2,\ldots,t$

**Corollary 5.1.2.** *Let* $a+ih$ $(i=-m,\ldots-2,-1,0,1,2,\ldots n)$ *are the distinct numbers spaced equally in the interval* $[p,q]$, $h\neq 0$, $t\in W$ *and* $f\in C^{n+m+t}[p,q]$, *then for some* $\xi_i\in\{a-mh,\ldots,a-h,a,a+h,\ldots a+nh\}$, $i=0,1,\ldots,t-1$

$$\frac{f^{(t)}(a)}{t!}=\frac{1}{h^t}\left(\sum_{i=1}^{n}A_i\left(\frac{\hat{a}_0}{i^t}+\frac{\hat{a}_1}{i^{t-1}}+\ldots+\frac{\hat{a}_{t-1}}{i}\right)f(a+ih)\right.$$
$$+\sum_{i=1}^{m}A_{-i}\left(\frac{(-1)^t\hat{a}_0}{i^t}+\frac{(-1)^{t-1}\hat{a}_1}{i^{t-1}}+\ldots-\frac{\hat{a}_{t-1}}{i}\right)f(a-ih)+\hat{a}_t f(a)\Bigg) \qquad , \ t=0,1,2,\ldots \quad (5.11)$$
$$+(-1)^n h^{n+m} n!m!\left(\hat{a}_0\frac{f^{(t+n+m)}(\xi_0)}{(t+n+m)!}+\frac{\hat{a}_1}{h}\frac{f^{(t+n+m-1)}(\xi_1)}{(t+n+m-1)!}+\ldots+\frac{\hat{a}_{t-1}}{h^{t-1}}\frac{f^{(n+m+1)}(\xi_{t-1})}{(n+m+1)!}\right)\Bigg)$$

Where $W^{(r)}=\sum_{i=1}^{m}\frac{(-1)^r A_{-i}}{i^r}+\sum_{i=1}^{n}\frac{A_i}{i^r}$, $r=1,2,\ldots t$

$A_{-i}=\frac{(-1)^{i-1}(i+1)(i+2)\ldots m.n!}{(m-i)!.(i+1)(i+2)\ldots(i+n)}$, $i=1,2,\ldots m$ and $A_i=\frac{(-1)^{i-1}m!.(i+1)(i+2)\ldots n}{(i+1)(i+2)\ldots(i+m)(n-i)!}$ , $i=1,2,\ldots n$

$\hat{a}_0=1$, $\hat{a}_k=-(W^{(k)}\hat{a}_0+W^{(k-1)}\hat{a}_1+W^{(k-2)}\hat{a}_2+\ldots+W^{(1)}\hat{a}_{k-1})$, $k=1,2,\ldots,t$

**Corollary 5.1.3.** *Let* $a\pm ih$ $(i=0,1,2,\ldots n)$ *are the distinct numbers spaced equally in the interval* $[p,q]$, $h\neq 0$, *and* $f\in C^{2n+t}[p,q]$, *then there exists a number* $\xi_i\in\{a-nh,\ldots,a-h,a,a+h,\ldots a+nh\}$ *with* $i=0,1,\ldots,t-1$ *then*

$$\frac{f^{(t)}(a)}{t!}=\frac{1}{h^t}\left(\sum_{i=1}^{n}A_i\left(f(a+ih)+(-1)^{\psi}f(a-ih)\right)\left(\frac{\tilde{a}_0}{i^t}+\frac{\tilde{a}_2}{i^{t-2}}+\ldots\frac{\tilde{a}_{t-2}}{i^2}\right)+\tilde{a}_t f(a)\right.$$
$$+h^{2n}(-1)^n n!n!\left(\tilde{a}_0\frac{f^{(t+2n)}(\xi_0)}{(t+2n)!}+\frac{\tilde{a}_1}{h}\frac{f^{(t+2n-1)}(\xi_1)}{(t+2n-1)!}+\ldots+\frac{\tilde{a}_{t-1}}{h^{t-1}}\frac{f^{(2n+1)}(\xi_{t-1})}{2n+1!}\right)\Bigg) \quad , \ t=0,1,2,\ldots \qquad (5.12)$$

Where $A_i=(-1)^{i-1}\frac{n(n-1)\ldots(n-i+1)}{(n+1)(n+2)\ldots(n+i)}$, $i=1,2,3,\ldots,n$

$V^{(k)}=\sum_{i=1}^{n}\frac{A_i}{i^k}$, $k=2,4,6,\ldots$ $\psi=\begin{cases}0, & k \text{ is even}\\ 1, & k \text{ is odd}\end{cases}$

and $\tilde{a}_0=1$, $\tilde{a}_k=-2(V^{(2)}\tilde{a}_{k-2}+V^{(4)}\tilde{a}_{k-4}+\ldots+V^{(k)}\tilde{a}_0)$, $k=2,4,6\ldots$

Corollary (5.1.2) is generalized formula for forward, backward and central difference formulas of numerical differentiation. If we put $n=0$, $m=0$ and $m=n$ then we get forward, backward and central difference formulas respectively. Also, we get more number of formulas for different values of *m* and *n*. For example, various five point formulas for second order differentiation is given below.

$$f^{(2)}(a)=\frac{1}{12h^2}\left(11f(a-4h)-56f(a-3h)+114f(a-2h)-104f(a-h)+35f(a)\right)+O(h^3) \qquad (5.13)$$

$$f^{(2)}(a)=\frac{1}{12h^2}\left(-f(a-3h)+4f(a-2h)+6f(a-h)-20f(a)+11f(a+h)\right)+O(h^3) \qquad (5.14)$$

$$f^{(2)}(a)=\frac{1}{12h^2}\left(-f(a-2h)+16f(a-h)-30f(a)+16f(a+h)-f(a+2h)\right)+O(h^4) \qquad (5.15)$$

$$f^{(2)}(a)=\frac{1}{12h^2}\big(11f(a-h)-20f(a)+6f(a+h)+4f(a+2h)-f(a+3h)\big)+O(h^3) \qquad (5.16)$$

$$f^{(2)}(a)=\frac{1}{12h^2}\big(35f(a)-104f(a+h)+114f(a+2h)-56f(a+3h)+11f(a+4h)\big)+O(h^3) \qquad (5.17)$$

Similarly, the higher order derivatives are easily found from the above formulas. It can be easily found that these formulas are the same as to those known corresponding numerical differentiation formulas based on interpolating polynomials (such as the Lagrange, Newton interpolating polynomials). But, the procedure of evaluating higher order derivative differs in different methods. However, the best method is chosen on the basis of computation time and storage space requirement. As mentioned before, the forms based on interpolating polynomials and other former formulas are complicated, whereas the new method here has some important advantages. It gives formulas that use given function values at the points directly and easily to calculate numerical approximations of arbitrary order at any sampling data for higher derivatives. This new method need less calculation burden, computing time and storage to estimate the derivatives than the other methods stated above. For example, suppose $x, x_0, x_1, \ldots, x_n$ are distinct numbers in the interval $[p,q]$, then the number of various arithmetic operations on the data, to evaluate $f^{(k)}(x)/k!$ using Lagrange interpolating polynomial, Theorem 2.2 in Ref [13] and Equation (5.5) in this paper are compared and given below. If *n* and *k* increases, then the total number of operations rapidly increases in both Lagrange interpolating polynomial and Theorem 2.2 in Ref[13]. But, the Equation (5.5) requires minimum number of operations then them

Table 9. Comparisons of Total number of different operations with former formulas.

| Operations | Based on Lagrange Polynomial, $k \geq 0$ | Theorem 2.2 In Ref[13], $k \geq 1$ | Equation (5.5), $k \geq 1$ |
|---|---|---|---|
| Additions | $(n+1)^n C_k - 1$ | $(n+1)^n C_{k-1} - 1$ | $n(2k+1)+k(k+1)/2$ |
| Subtractions | $(n+1)\big({}^nC_k(n-k)+n\big)$ | $(n+1)\big({}^nC_{k-1}(n-k+1)+n+2\big)$ | $(n+1)(2n+k^2+k)$ |
| Multiplications | $(n+1)\big({}^nC_k(n-k-1)+n\big)$ | $(n+1)\big({}^nC_{k-1}(n-k)+n\big)$ | $2n(n+1)+(n+1)k(k-1)+k(k+1)/2$ |
| Division | $n+1$ | $2(n+1)$ | $(n+1)(2k+1)$ |

Further, in section 4, we have studied k^th^ order differentiation formula in terms of linear combination of all possible k^th^ order divided differences and simple coefficients. It is a direct way to evaluate higher order differentiation, but the problem with this method is, it requires all combination of divided differences. However, the coefficients of divided differences are very simple polynomials. For example, if we want to find third order differentiation formula we need to evaluate all possible 3^rd^ order divided differences if the functional value is not known. But, this method is very easy to generalize to arbitrary number of data for evenly and unevenly spaced points.

**Algorithm 5.1.** For the unevenly spaced points $x_0, x_1, x_2, \ldots, x_n$ and known the functional values $f(x_i)$ at $x_i$, $i = 0,1,2,\ldots,n$ then the steps to use ($n$+1) point formula to estimate $t^{th}$ derivative of *f*(*x*) at *x* are

Step1: For $i = 0$ *to* $n$ do $L_i(x)=\prod_{j=0,i\neq j}^{n}\frac{x-x_j}{x_i-x_j}$

Step 2: For $k=1$ *to* $t$ do $\rho_k=\sum_{i=0}^{n}\frac{L_i(x)}{(x_i-x)^k}$

Step3: $a_0 = 1$, For $k=1$ *to* $t$ do $a_k = -(\rho_k a_0 + \rho_{k-1}a_1 + \rho_{k-2}a_2 + \ldots + \rho_1 a_{k-1})$

Step 4: Use (5.5) to find $t^{th}$ derivative at '*x*'.

**Algorithm 5.2.** For the evenly spaced points with interval $h$ and known the functional values $f(x_i)$ at $x_i = a + ih$, $i = -m,\ldots,-2,-1,0,1,2,\ldots,n$ then the steps to estimate $t^{th}$ derivative of $f(a)$ at $a$, are

Step 1: For $i = 1$ *to* $m$ do $A_{-i} = \dfrac{(-1)^{i-1}(i+1)(i+2)\ldots m.n!}{(m-i)!(i+1)(i+2)\ldots(i+n)}$

Step 2: For $i = 1$ *to* $n$ do $A_i = \dfrac{(-1)^{i-1}m!(i+1)(i+2)\ldots n}{(i+1)(i+2)\ldots(i+m).(n-i)!}$

Step3: For $r = 1$ *to* $t$ do $W^{(r)} = \sum_{i=1}^{m}\dfrac{(-1)^r A_{-i}}{i^r} + \sum_{i=1}^{n}\dfrac{A_i}{i^r}$

Step 4: $\hat{a}_0 = 1$, For $k = 1$ *to* $t$ do $\hat{a}_k = -(W^{(k)}\hat{a}_0 + W^{(k-1)}\hat{a}_1 + W^{(k-2)}\hat{a}_2 + \ldots + W^{(1)}\hat{a}_{k-1})$

Step 5: Use (5.11) to find $t^{th}$ derivative at '$a$'.

## 6. Formulas for Integration

**Theorem 6.1.** *Let* $x_0, x_1, x_2, \ldots, x_n$ *are the distinct numbers in the interval* $[p,q]$, $x, x+h \in [p,q]$, $h \neq 0$ *and if* $f \in C^{2n+1}[p,q]$, *with*

$$\int_x^{x+h} f(x)dx = \sum_{i=0}^{n} f(x_i)L_i(x)\left(\gamma_0 + \frac{\gamma_1}{(x_i - x)} + \frac{\gamma_2}{(x_i - x)^2} + \ldots + \frac{\gamma_n}{(x_i - x)^n}\right)$$
$$+ \prod_{i=0}^{n}(x - x_i)\left(\frac{f^{(n+1)}(\xi_n)}{(n+1)!}\gamma_0 + \frac{f^{(n+2)}(\xi_{n-1})}{n+2!}\gamma_1 + \ldots + \frac{f^{(2n+1)}(\xi_0)}{(2n+1)!}\gamma_n\right) + O(h^{n+2}) \tag{6.1}$$

Where $L_i(x_0) = \prod_{\substack{j=0 \\ i \neq j}}^{n}\dfrac{x - x_j}{x_i - x_j}$, $i = 0,1,2,\ldots,n$, $\rho_k = \sum_{i=0}^{n}\dfrac{L_i(x)}{(x_i - x)^k}$, $k = 0,1,2,\ldots n$

$a_0 = 1$, $a_k = -(\rho_k a_0 + \rho_{k-1}a_1 + \rho_{k-2}a_2 + \ldots + \rho_1 a_{k-1})$, $k = 0,1,2,\ldots n$

$\gamma_k = \dfrac{h^{k+1}}{k+1} + a_1\dfrac{h^{k+2}}{k+2} + \cdots + a_{n-k}\dfrac{h^{n+1}}{n+1}$, $k = 0,1,2,\ldots n$ and $\xi_i \in \{x, x_0, x_1, \ldots x_n\}$, $i = 0,1,2,\ldots,n$

**Proof.**

Using Taylor series on integration

$$\int_x^{x+h} f(x)dx = f(x)h + f'(x)\frac{h^2}{2!} + f''(x)\frac{h^3}{3!} + \ldots + f^{(n)}(x)\frac{h^{n+1}}{(n+1)!} + O(h^{n+2}) \tag{6.2}$$

Using (5.1) in equation (6.2),

$$\int_x^{x+h} f(x)dx = f(x)h + \big(M_1(x) + f(x)a_1\big)\frac{h^2}{2} + \big(M_2(x) + M_1(x)a_1 + f(x)a_2\big)\frac{h^3}{3} +$$
$$\ldots + \big(M_n(x) + M_{n-1}(x)a_1 + M_{n-2}(x)a_2 + \ldots + M_1(x)a_{n-1} + f(x)a_n\big)\frac{h^{n+1}}{n+1} + O(h^{n+2}) \tag{6.3}$$

Rearranging the above equation,

$$
\begin{aligned}
&= f(x)\left(h + a_1\frac{h^2}{2} + a_2\frac{h^3}{3} + \cdots + a_n\frac{h^{n+1}}{n+1}\right) + M_1(x)\left(\frac{h^2}{2} + a_1\frac{h^3}{3} + \cdots + a_{n-1}\frac{h^{n+1}}{n+1}\right) \\
&\quad + M_2(x)\left(\frac{h^3}{3} + a_1\frac{h^4}{4} + \cdots + a_{n-2}\frac{h^{n+1}}{n+1}\right) + \ldots + M_n(x)\frac{h^{n+1}}{n+1} + O(h^{n+2})
\end{aligned}
\tag{6.4}
$$

$$= f(x)\gamma_0 + M_1(x)\gamma_1 + M_2(x)\gamma_2 + \ldots + M_n(x)\gamma_n + O(h^{n+2})$$

But $M_k(x) = \sum_{i=0}^{n} \frac{f(x_i)}{(x_i - x)^k} L_i(x) + \frac{f^{(k+n+1)}(\xi_{n-k})}{(k+n+1)!} \prod_{i=0}^{n}(x - x_i)$, $\xi_k \in \{x, x_0, x_1, \ldots x_n\}$, $k = 0,1,2,\ldots,n$

Using Lagrange interpolation formula and substituting all *M*'s in (6.4) and after simplification, we get

$$
\begin{aligned}
&= \sum_{i=0}^{n} f(x_i)L_i(x)\gamma_0 + \sum_{i=0}^{n} \frac{f(x_i)}{(x_i - x)} L_i(x)\gamma_1 + \ldots + \sum_{i=0}^{n} \frac{f(x_i)}{(x_i - x)^n} L_i(x)\gamma_n \\
&\quad + \prod_{i=0}^{n}(x - x_i)\left(\frac{f^{(n+1)}(\xi_n)}{(n+1)!}\gamma_0 + \frac{f^{(n+2)}(\xi_{n-1})}{n+2!}\gamma_1 + \ldots + \frac{f^{(2n+1)}(\xi_0)}{(2n+1)!}\gamma_n\right) + O(h^{n+2})
\end{aligned}
$$

$$
\begin{aligned}
&= \sum_{i=0}^{n} f(x_i)L_i(x)\left(\gamma_0 + \frac{\gamma_1}{(x_i - x)} + \ldots + \frac{\gamma_n}{(x_i - x)^n}\right) \\
&\quad + \prod_{i=0}^{n}(x - x_i)\left(\frac{f^{(n+1)}(\xi_0)}{(n+1)!}\gamma_0 + \frac{f^{(n+2)}(\xi_1)}{n+2!}\gamma_1 + \ldots + \frac{f^{(2n+1)}(\xi_n)}{(2n+1)!}\gamma_n\right) + O(h^{n+2})
\end{aligned}
\tag{6.5}
$$

Hence the proof.

**Corollary 6.1.1.** *Let* $a, a+h, a+2h, \ldots, a+nh$ *are distinct numbers in the interval* $[p,q]$, *spaced equally with* $h \neq 0$ *and if* $f \in C^{2n}[p,q]$, *with*

$$\int_a^{a+nh} f(x)dx = h\left[f(a)\xi_0 + \sum_{i=1}^{n}(-1)^{i-1}\, {}^nC_i\, f(a+ih)\left(\frac{\xi_1}{i} + \frac{\xi_2}{i^2} + \ldots + \frac{\xi_n}{i^n}\right)\right] + O(h^{n+2}) \tag{6.6}$$

Where $U^{(k)} = \sum_{i=1}^{n} \frac{(-1)^{i-1}}{i^k}\, {}^nC_i$, $k = 1,2,\ldots,n$

$\hat{a}_0 = 1$, $\hat{a}_k = -(U^{(k)}\hat{a}_0 + U^{(k-1)}\hat{a}_1 + U^{(k-2)}\hat{a}_2 + \ldots + U^{(1)}\hat{a}_{k-1})$, $k = 1,2,\ldots,n$ and

$\xi_k = \frac{n^{k+1}}{k+1} + \hat{a}_1\frac{n^{k+2}}{k+2} + \cdots + \hat{a}_{n-k}\frac{n^{n+1}}{n+1}$ , $k = 0,1,2,\ldots,n$

**Proof.**

Using Taylor series on integration

$$\int_a^{a+nh} f(x)dx = f(a)nh + f'(a)\frac{(nh)^2}{2!} + f''(a)\frac{(nh)^3}{3!} + \ldots + f^{(n)}(a)\frac{(nh)^{n+1}}{(n+1)!} + O(h^{n+2}) \tag{6.7}$$

Using equation (5.2) then

$$\frac{f^{(t)}(x)}{t!} = \frac{1}{h^t}\left(\hat{M}_t(a) + \hat{M}_{t-1}(a)\hat{a}_1 + \ldots + \hat{M}_1(a)\hat{a}_{t-1} + f(a)\hat{a}_t\right),\ t = 1,2,\ldots,n \tag{6.8}$$

Where $\hat{M}_k(a) = \sum_{i=1}^{n} \frac{f(a+ih)}{i^k}(-1)^{i-1}\,{}^nC_i + (-1)^n h^{n+k}\frac{f^{(k+n)}(\xi_{n-k})}{(k+n)!}$, $k = 1,2,\ldots,n$

Using (6.8) in (6.7)

$$\begin{aligned}\int_a^{a+nh} f(x)dx &= f(a)nh + \left(\hat{M}_1(a) + \hat{a}_1\right)\frac{n^2h}{2} + \left(\hat{M}_2(a) + \hat{M}_1(a)\hat{a}_1 + f(a)\hat{a}_2\right)\frac{hn^3}{3} + \\ &\ldots + \left(\hat{M}_n(a) + \hat{M}_{n-1}(a)\hat{a}_1 + \hat{M}_{n-2}(a)\hat{a}_2 + \ldots + \hat{M}_1(a)\hat{a}_{n-1} + f(a)\hat{a}_n\right)\frac{hn^{n+1}}{n+1} + O(h^{n+2})\end{aligned} \tag{6.9}$$

$$\begin{aligned}&= h\Bigg[f(a)\left(n + \hat{a}_1\frac{n^2}{2} + \hat{a}_2\frac{n^3}{3} + \cdots + \hat{a}_n\frac{n^{n+1}}{n+1}\right) + \hat{M}_1(a)\left(\frac{n^2}{2} + \hat{a}_1\frac{n^3}{3} + \cdots + \hat{a}_{n-1}\frac{n^{n+1}}{n+1}\right) \\ &\quad + \hat{M}_2(a)\left(\frac{n^3}{3} + \hat{a}_1\frac{n^4}{4} + \cdots + \hat{a}_{n-2}\frac{n^{n+1}}{n+1}\right) + \ldots + \hat{M}_n(a)\frac{n^{n+1}}{n+1}\Bigg] + O(h^{n+2})\end{aligned} \tag{6.10}$$

Let $\xi_k = \frac{n^{k+1}}{k+1} + \hat{a}_1\frac{n^{k+2}}{k+2} + \cdots + \hat{a}_{n-k}\frac{n^{n+1}}{n+1}$ , $k = 0,1,2,\ldots,n$ then,

$$= h\left(f(a)\xi_0 + \hat{M}_1(a)\xi_1 + \hat{M}_2(a)\xi_2 + \ldots + \hat{M}_n(a)\xi_n\right) + O(h^{n+2})$$

Substituting all $\hat{M}$ 's values we obtain (6.6).

**Corollary 6.1.2.** *Let* $a \pm ih$ $(i = 0,1,2,\ldots n)$ *are the distinct numbers spaced equally in the interval* $[p,q]$, $h \neq 0$ *and* $f \in C^{3n}[p,q]$

$$\int_{a-nh}^{a+nh} f(x)dx = 2h\left[f(a)\xi_0 + \sum_{i=1}^{n}\left[f(a+ih) + f(a-ih)\right]A_i\left(\frac{\xi_2}{i^2} + \frac{\xi_4}{i^4}\ldots + \frac{\xi_{2n}}{i^{2n}}\right)\right] + O(h^{2n+2)} \tag{6.11}$$

Where $A_i = (-1)^{i-1}\frac{n(n-1)(n-2)\ldots(n-i+1)}{(n+1)(n+2)(n+3)\ldots(n+i)}$, $V^{(2k)} = \sum_{i=1}^{n}\frac{A_i}{i^{2k}}$ , $k = 1,2,\ldots,n$

$\tilde{a}_0 = 1$, $\tilde{a}_k = -2(V^{(2)}\tilde{a}_{k-2} + V^{(4)}\tilde{a}_{k-4} + \ldots + V^{(k-2)}\tilde{a}_2 + V^{(k)}\tilde{a}_0)$, $k = 2,4,6,\ldots 2n$ and

$\xi_{2k} = \frac{n^{2k+1}\tilde{a}_{2k}}{2k+1} + \frac{n^{2k+3}\tilde{a}_{2k+2}}{2k+3} + \ldots + \frac{n^{2n+1}\tilde{a}_{2n-2k}}{2n+1}$ , $k = 0,1,2,\ldots,n$

Corollary 6.1.1 is general formula of trapezoidal, Simpson and Newton cotes integration formulas. If $n = 2,6$ we have the two approximations for numerical integration

$$\int_{a}^{a+2h} f(x)dx \approx \frac{h}{3}\left(f(a) + 4f(a+h) + f(a+2h)\right) \tag{6.12}$$

$$\int_{a}^{a+6h} f(x)dx \approx \frac{h}{140}\left(41f(a) + 216f(a+h) + 27f(a+2h) + 272f(a+3h) + 27f(a+4h) + 216f(a+5h) + 41f(a+6h)\right) \tag{6.13}$$

**Algorithm 6.1.** For the unevenly spaced points $x_0, x_1, x_2, \dots, x_n$ and known the functional values $f(x_i)$ at $x_i$, $i = 0,1,2,\dots,n$ then the steps to estimate $\int_{x}^{x+h} f(x)dx$ are,

Step1: For $i = 0$ *to* $n$ do $L_i(x) = \prod_{j=0, i\neq j}^{n} \frac{x - x_j}{x_i - x_j}$

Step 2: For $k = 1$ *to* $n$ do $\rho_k = \sum_{i=0}^{n} \frac{L_i(x)}{(x_i - x)^k}$,

Step3: $a_0 = 1$, For $k = 1$ *to* $n$ do $a_k = -(\rho_k a_0 + \rho_{k-1} a_1 + \rho_{k-2} a_2 + \dots + \rho_1 a_{k-1})$

Step 4: For $k = 0$ *to* $n$ do $\gamma_k = \frac{h^{k+1}}{k+1} + a_1 \frac{h^{k+2}}{k+2} + \cdots + a_{n-k} \frac{h^{n+1}}{n+1}$

Step 5: Use (6.1), to find numerically $\int_{x}^{x+h} f(x)dx$

**Algorithm 6.2.** For the evenly spaced points with interval *h,* known the functional values $f(x_i)$ at $x_i = a + ih$, $i = 1,2,\dots,n$ then the steps to estimate $\int_{a}^{a+nh} f(x)dx$ are,

Step 1: For $k = 1$ *to* $n$ do $U^{(k)} = \sum_{i=1}^{n} \frac{(-1)^{i-1}}{i^k} {}^{n}C_i$

Step 2: $\hat{a}_0 = 1$, For $k = 1$ *to* $n$ do $\hat{a}_k = -(U^{(k)} \hat{a}_0 + U^{(k-1)} \hat{a}_1 + U^{(k-2)} \hat{a}_2 + \dots + U^{(1)} \hat{a}_{k-1})$

Step3: For $k = 0$ *to* $n$ do $\xi_k = \frac{n^{k+1}}{k+1} + \hat{a}_1 \frac{n^{k+2}}{k+2} + \cdots + \hat{a}_{n-k} \frac{n^{n+1}}{n+1}$

Step 4: Use (6.6), to find numerically $\int_{a}^{a+nh} f(x)dx$

## 7. Conclusion

By introducing a new formula of divided difference and new schemes of divided difference tables, we have studied interpolation, numerical differentiation and integration formulas with arbitrary order accuracy for evenly and unevenly spaced data. First, we study the new interpolation formula which generalizes both Newton's and Lagrange's interpolation formulas and its various attractive forms. Further, we study other new interpolation formulas based on the differences and divided differences. Also comparisons of former interpolation formulas (Newton's, Gauss's, Stirling, Bessel's, etc.,) with the new modified formulas are shown that the new formulas are very efficient and posses good accuracy for evaluating functional values between given data. Comparison table (Table 8) in section 3 shows the new interpolation formula is superior than Newton's and Lagrange's interpolation formulas.

Second, we study two new types of numerical differentiation formulas. In the first type, the higher order derivatives are evaluated recursively by using the new recursive formulas for both equally and unequally spaced data in section 4. The new recursive formula is more general formula for both interpolation and numerical differentiation and given with various remainder terms. Using this, we find higher order derivatives in terms of functional values at given nodes in section 5. Also, in the same section, the total number of various operations is compared between former formulas and the new formulas and this show the new formulas require very less computation time and storage space than the former formulas. In the second type, we study another way of evaluating numerical differentiation using linear combination of divided differences for evenly and unevenly spaced data in various cases. It requires more number of divided differences. However it is very easy to program in computer. Third, the numerical integration formulas are studied for both equally and unequally spaced data to any arbitrary order accuracy. Basic computer algorithms are given for new formulas. Thus, the three major problems in numerical analysis are solved by new formulas and methods through new divided difference formula and new divided difference tables in this paper.

**List of References**